\documentclass[12pt,a4paper]{article}

\usepackage{amsmath}
\usepackage{amssymb}
\usepackage{amsthm}
\usepackage{bbm}
\usepackage{enumerate}
\usepackage{latexsym}
\usepackage{mathdots}
\usepackage{geometry}

\usepackage{hyperref}   

\theoremstyle{definition}
\newtheorem{mydefn}{Definition}[section]
\newtheorem{remark}{Remark}[section]

\theoremstyle{plain}
\newtheorem{theorem}{Theorem}[section]
\newtheorem{lemma}{Lemma}[section]
\newtheorem{corollary}{Corollary}[section]
\newtheorem{proposition}{Proposition}[section]
\newtheorem{conjecture}{Conjecture}[section]

\newcommand{\tpi}{2\pi i}

\newcommand{\Z}{\mathbb{Z}}

\newcommand{\R}{\mathbb{R}}
\newcommand{\C}{\mathbb{C}}
\newcommand{\HH}{\mathbb{H}}

\DeclareMathOperator{\SL}{SL}
\DeclareMathOperator{\Sp}{Sp}
\DeclareMathOperator{\End}{End}
\DeclareMathOperator{\wt}{wt}
\DeclareMathOperator{\Tr}{Tr}

\newcommand{\omt}{\widetilde{\omega}}
\newcommand{\vac}{\mathbf{1}}

\newcommand{\abar}{\overline{a}}
\newcommand{\sups}[1]{\,^{#1}\hspace{-0.5 mm}}
\newcommand{\supsin}[1]{\,^{#1}\hspace{-1 mm}}

\begin{document}


\title{Genus Two Zhu Theory for \\Vertex Operator Algebras}
\author{
Thomas Gilroy\thanks{Supported by an Irish Research Council Government of Ireland Studentship}
\\
School of Mathematics and Statistics,\\
University College Dublin,\\
Dublin 4, Ireland
\and
Michael P. Tuite 
\\
School of Mathematics, Statistics and Applied Mathematics, \\
National University of Ireland Galway, \\
University Road, Galway, Ireland
}

\maketitle


\begin{abstract}
\noindent
We consider correlation functions for a vertex operator algebra on a genus two Riemann surface formed by sewing two tori together. 
We describe a generalisation of genus one Zhu recursion where we express an arbitrary genus two $n$--point correlation function in terms of $(n-1)$--point functions. 
We consider several applications  including the correlation functions for the Heisenberg vertex operator algebra and its modules,   Virasoro correlation functions and genus two Ward identities. 
We  derive novel differential equations in terms of a differential operator on the genus two Siegel upper half plane for holomorphic $1$--differentials, the normalised bidifferential of the second kind, the projective connection and the Heisenberg partition function. We  prove that the holomorphic mapping from the sewing parameter domain to the Siegel upper half plane is injective but not surjective. We also demonstrate that genus two differential equations arising from Virasoro singular vectors have holomorphic coefficients.
\end{abstract}


\section{Introduction}
\noindent
The connection between Vertex Operator Algebras (VOAs) and elliptic functions and modular forms has been a fundamental aspect of the theory since its inception  in the work of Borcherds \cite{B} and Frenkel, Lepowsky and Meurmann \cite{FLM}. This phenomenon is manifested through  $n$--point correlation trace functions. Zhu recursion expresses a genus one correlation $n$-point  function in terms of $(n-1)$--point functions using a formal recursive identity involving Weierstrass elliptic functions \cite{Z}. Zhu recursion implies genus one trace functions satisfy modular differential equations if the VOA is $C_2$--cofinite  \cite{Z}. Such differential equations imply convergence and modular properties for trace functions. Zhu recursion is also an important calculation tool e.g. all correlation functions for the Heisenberg VOA and its modules can be computed exactly \cite{MT1}.

The expression of genus one correlation functions in terms of  elliptic and modular functions is also fundamental in conformal field theory \cite{BPZ, DMS}. 
The importance of extending these ideas to a general genus Riemann surface was also recognised early on in physics e.g. \cite{Kn, FS}.
There are also natural mathematical reasons to extend this connection to Riemann surfaces of higher genus. In particular, we would like to understand how the elliptic functions and modular forms of genus one Zhu theory generalise at higher genus and develop a scheme to study the convergence and modular properties of genus two partition functions via genus two differential equations. 

In recent work,
correlation functions for VOAs and super-VOAs on a genus two Riemann surface have been defined, and in some cases computed \cite{T,MT3,MT4,TZ2,TZ3}, based on explicit sewing procedures \cite{Y,MT2, TZ1}. The present paper deals with correlation functions for VOAs on a genus two Riemann genus two surface formed by sewing two tori together. In particular, we describe a new formal Zhu recursion formula expressing any genus two $n$--point function in terms of $(n-1)$--point functions. We consider a number of applications paralleling the genus one case.

In Sect.~\ref{sec: Riemannreview} we review the relevant parts of the complex analytic theory of Riemann surfaces. We begin with definitions of elliptic functions and modular forms. We then discuss aspects of general genus $g$ Riemann surfaces including a modular invariant differential operator on the Siegel upper half plane and introduce a general genus analogue of the Serre derivative exploited in later sections. We conclude  with a discussion on the construction of a genus two Riemann surface by sewing together two tori.

In Sect.~\ref{sec: VOA12review} we review relevant definitions and results regarding vertex operator algebras and we review Zhu recursion for genus one correlation functions. The sewing procedure of Sect.~\ref{sec: Riemannreview} informs the definition of genus two VOA correlation functions in terms of infinite formal sums of appropriate genus one correlation functions \cite{T,MT3,MT4}.

Sect.~\ref{sec:Zhured} contains the new genus two formal Zhu recursion identity relating genus two $n$--point correlation functions to $(n-1)$--point functions. This  is achieved by applying genus one Zhu reduction to the genus one component parts  leading to a system of recursive identities which are  solved to obtain genus two Zhu recursion. 
The expansion in  $(n-1)$--point functions uses a family of new generalised genus two Weierstrass formal functions analogous to the elliptic Weierstrass functions appearing in genus one Zhu recursion. 
Unlike the genus one case, the generalised  Weierstrass functions depend on the conformal weight $N$ of the vector being reduced but are otherwise universal. 
Another novel feature of genus two Zhu recursion is that a genus two $1$-point function of a weight $N$ vector is expressed in terms of 2 coefficient functions for $N=1$ and $2N-1$ coefficient functions for $N\ge 2$. This agrees with the dimension of the space of genus two holomorphic $N$--differentials according to the Riemann-Roch theorem.

In Sect.~\ref{sec:Zhuwt1} we consider genus two Zhu recursion in the case where we reduce on a vector of weight $N=1$. We prove the holomorphy of the coefficient functions in terms of holomorphic 1--differentials and the generalised Weierstrass functions in terms of the normalised bidifferential form of the second kind. We calculate the genus two Heisenberg $n$-point correlation functions for a pair of Heisenberg modules. This agrees with results of \cite{MT3} obtained by  combinatorial methods. 

Sect.~\ref{sec:Zhuwt2} we consider genus two Zhu recursion for a vector of weight $N=2$. 
We prove that the coefficient terms for the genus two $1$-point function are holomorphic 2--differentials. 
We show that the genus two Virasoro $1$--point function is given by a certain derivative of the partition function with respect to the sewing parameters.
We also derive the genus two $n$--point correlation functions for $n$ Virasoro vectors and a genus two Ward Identity where the new derivative again plays a role. 

In Sect.~\ref{sec:GDE} we relate the  differential operator of Sect.~\ref{sec:Zhuwt2} to the modular invariant differential operator of Sect.~\ref{sec: Riemannreview}. We obtain a closed holomorphic formula for the $N=2$ generalised Weierstrass function. This completes the proof of the holomorphy of all coefficient terms appearing in the genus two Ward identities and Virasoro $n$-point functions of Sect.~\ref{sec:Zhuwt2}. As an important consequence, this implies that differential equations arising from Virasoso singular vectors therefore have holomorphic coefficients. These developments are further explored in \cite{GT1,GT2}.   
We also prove that the holomorphic map from the sewing domain to the Siegel upper half plane is injective but not surjective. Finally, we describe novel holomorphic differential equations for  $1$--differentials, the normalised $2$--bidifferential, the projective connection and the genus two Heisenberg partition function.


\section{Review of Riemann Surfaces}
\label{sec: Riemannreview}

We begin with some basic notations and definitions that will be used throughout the paper. $\Z,\R$ and $\C$ denote the integers, reals and complex numbers respectively. 
$
\HH = \{ \tau \in \C \vert \Im(\tau) > 0\}$
 is the complex upper-half plane. 
We use the conventions that $q = e^{2\pi i \tau}$ for $\tau\in \HH$ and $q_x = e^x$. 
For derivatives we use the shorthand notation $\partial_x  = \frac{\partial}{\partial x}$.

\subsection{Elliptic functions and modular forms}
\label{subsect_Elliptic}
We define some elliptic functions and modular forms \cite{Se,La}.
\begin{mydefn}\label{def:Eisenstein}
The \emph{Eisenstein series} for an integer $k\geq 2$ is given by
\begin{align*}
E_k(\tau) = E_k(q) = \begin{cases}
0, & \mbox{ for }k \mbox{ odd}, \\
\displaystyle -\frac{B_{k}}{k!}+\frac{2}{(k-1)!}\sum_{n\geq 1}\sigma
_{k-1}(n)q^{n}, & \mbox{ for }k \mbox{ even}.
\end{cases}
\end{align*}
where $\tau\in \HH$ ($q=e^{2\pi i \tau}$), $\sigma_{k-1}(n) = \sum_{d\vert n} d^{k-1}$ and $k^{\mathrm{th}}$ Bernoulli number  $B_k$.
\end{mydefn}
\noindent
If $k\geq 4$ then $E_k(\tau)$ is a  modular form of weight $k$ on $\SL(2,\Z)$, while $E_2(\tau)$ is a quasi-modular form. We also define elliptic functions $z\in \C$
\begin{mydefn}\label{def:Pfunctions}
For integer $k\ge 1$
\begin{align*}
P_{k}(z,\tau)=\frac{(-1)^{k-1}}{(k-1)!}\partial_z^{k-1}P_{1}(z,\tau),\quad P_{1}(z,\tau)=\frac{1}{z}-\sum_{k\geq 2}E_{k}(\tau )z^{k-1}.
\end{align*}
\end{mydefn}
\noindent
In particular $P_{2}(z, \tau) =\wp (z,\tau)+E_{2}(\tau )$ for 
Weierstrass function $\wp (z,\tau)$ with periods $2\pi i$ and $2\pi i\tau$.
$P_{1}(z, \tau)$ is related to the quasi--periodic Weierstrass $\sigma$--function  with 
$P_{1}(z+2\pi i\tau, \tau)=P_{1}(z, \tau)-1$.


\subsection{Genus $g$ Riemann surfaces}
We review some relevant aspects of Riemann surface theory e.g. \cite{FK,F,Mu}. 
Consider a compact Riemann surface $\mathcal{S}^{(g)}$ of
genus $g$ with canonical homology cycle basis $\alpha^{i},\beta^i$ for $i=1,\ldots ,g$.
 There
exists $g$ holomorphic 1--differentials $\nu_{i}$   normalized by 
\begin{align}\label{nu}
\oint_{\alpha^{i}}\nu_{j}=2\pi i\delta_{ij}.
\end{align} 
These differentials can also be defined via the unique holomorphic bidifferential $(1,1)$--form $\omega(x,y)$ for $x\neq y$, known as the \emph{normalised bidifferential of the second kind}. It is defined by the following properties 
\begin{align}\label{omega}
\omega(x,y)=\left(\frac{1}{(x-y)^{2}}+\text{regular terms}\right)dxdy,
\end{align}
for any local coordinates $x,y$, with normalization 
\begin{align}\label{omeganorm}
\int_{\alpha^{i}}\omega(x,\cdot )=0,
\end{align}
for $i=1,\ldots ,g$.  Using the Riemann bilinear relations, one finds that 
\begin{align}\label{nui}
\nu_{i}(x)=\oint_{\beta^{i}}\omega(x,\cdot ),  
\end{align}
with  normalisation \eqref{nu}.
We also define the  \emph{period matrix} $\Omega$ by 
\begin{equation}
\Omega_{ij}=\frac{1}{2\pi i}\oint_{\beta^{i}}\nu_{j}, 
 \label{period}
\end{equation}%
for $i,j=1,\ldots ,g$ where $\Omega_{ij}=\Omega_{ji}$ and $\Im(\Omega)>0$ i.e. $\Omega \in \HH_{g}$, the Siegel upper half plane. 
Under a change of homology basis 
\begin{align}
\begin{bmatrix}
\widetilde{\boldsymbol{\beta}}\\
\widetilde{\boldsymbol{\alpha}}
\end{bmatrix}
=
\begin{bmatrix}
A & B\\
C & D 
\end{bmatrix}
\begin{bmatrix}
{\boldsymbol{\beta}}\\
{\boldsymbol{\alpha}}
\end{bmatrix},
\label{eq:hommod}
\end{align}
for $\left[
\begin{smallmatrix}
A & B\\
C & D 
\end{smallmatrix}
\right]\in\Sp(2g,\Z)$, the row vector 
$\boldsymbol{\nu}(x)=\left[\nu_1(x),\ldots ,\nu_g(x)\right] $ transforms to
\begin{align}
\widetilde{\boldsymbol{\nu}}(x)=\boldsymbol{\nu}M^{-1},
\label{eq:numod}
\end{align}
for $M=C\Omega+D$ where
the period matrix transforms to
\begin{align}
\widetilde{\Omega}=(A\Omega +B)(C\Omega+D)^{-1},
\label{eq:Omegamod}
\end{align}
and $\omega$ transforms to 
\begin{align}
\widetilde{\omega}(x,y)=&\omega(x,y)
-\frac{1}{2}\sum_{1\le i\le j\le g}\left(\nu_i(x)\nu_j(y)+\nu_j(x)\nu_i(y)\right)
\frac{\partial \log\det M}{\partial \Omega_{ij}}.
\label{eq:ommod}
\end{align} 
We also define the genus $g$ \emph{projective connection} $s(x)$  by
\begin{align}
s(x) = 6 \lim_{y\rightarrow x} \left( \omega(x,y) - \frac{dx dy}{(x-y)^2} \right).
\label{eq:projcon}
\end{align}
Under a conformal map $x\rightarrow \phi(x)$ one finds\footnote{The conventional factor of $6$  in the definition of $s(x)$ is introduced to simplify \eqref{eq:projSchwarz}.}  
\begin{align}
s(x)=s(\phi(x)) +\{\phi(x),x\}dx^2,
\label{eq:projSchwarz}
\end{align}
for the Schwarzian derivative $\{\phi(x),x\}=\frac{\phi'''(x)}{\phi'(x)}-\frac{3}{2}\left(\frac{\phi''(x)}{\phi'(x)} \right)^2$.
$s(x)$ is called a \emph{projective  form} since it transforms as a  $2$--differential under a M\"obius map $\phi(x)=\frac{ax+b}{cx+d}$ for which $\{\phi(x),x\}=0$.
From \eqref{eq:ommod}, 
$s(x)$ transforms under the modular group $\Sp(2g,\Z)$ to
\begin{align}
\widetilde{s} (x)=s(x)
-6\,\nabla_x \log\det M,
\label{eq:projcmod}
\end{align}
where $\nabla_x$ is the differential operator
\begin{align}
\nabla_x=
\sum_{1\le i \le j\le g}\nu_i(x)\nu_j(x)
\frac{\partial}{\partial \Omega_{ij}}.
\label{eq:nabla}
\end{align}
The subscript  indicates the dependence on $x\in \mathcal{S}^{(g)}$. The operator $\nabla_x$  will play an important role later on in this paper for genus $g=2$. 
\begin{proposition}\label{prop:nablaxmod}
\begin{align}
\frac{\partial\widetilde{\Omega}_{ab}}{\partial \Omega_{ij}}=
\begin{cases}
N_{ia}N_{ib} & i=j,\\
N_{ia}N_{jb} + N_{ib}N_{ja}& i\neq j,
\end{cases}
\label{eq:delOmtilde}
\end{align}
where $N=M^{-1}=(C\Omega+D)^{-1}$.  $\nabla_x$ is $\Sp(2g,\Z)$ invariant.
\end{proposition}
\noindent \textbf{Proof.} 
Using the $\Sp(2g,\Z)$ relations $A^TC=C^TA$  and $A^TD-C^TB=I_{g}$ we  note that   
\begin{align*}
N=A^T-C^T\widetilde{\Omega}.
\end{align*}
Consider 
\begin{align*}
\sum_{b=1}^{g}\frac{\partial\widetilde{\Omega}_{ab}}{\partial \Omega_{ij}} M_{bd}=
\frac{\partial(A\Omega+B)_{ad}}{\partial \Omega_{ij}}
-\sum_{b=1}^{g}\widetilde{\Omega}_{ab}\frac{\partial M_{bd}}{\partial \Omega_{ij}}.
\end{align*}
But
\begin{align*}
\frac{\partial}{\partial \Omega_{ij}}(A\Omega+B)_{ad}=
\begin{cases}
A_{ai}\delta_{di},& i=j,\\
A_{ai}\delta_{dj}+A_{aj}\delta_{di},& i\neq j,
\end{cases}
\end{align*}
with a similar formula for $\frac{\partial}{\partial \Omega_{ij}}M_{bd}$. This implies that for $i=j$
\begin{align*}
\sum_{a,b=1}^{g}\frac{\partial\widetilde{\Omega}_{ab}}{\partial \Omega_{ii}}M_{ac} M_{bd}
&=
\left((A^T-C^T\widetilde{\Omega})M\right)_{ic}\delta_{id}
=
\delta_{ic}\delta_{id}.
\end{align*}
Similarly, for $i\neq j$ we have
\begin{align*}
\sum_{a,b=1}^{g}\frac{\partial\widetilde{\Omega}_{ab}}{\partial \Omega_{ij}}M_{ac} M_{bd}
=
\delta_{ic}\delta_{jd}+\delta_{jc}\delta_{id}.
\end{align*}
Thus \eqref{eq:delOmtilde} follows. This immediately implies that  $\nabla_x \widetilde{\Omega}_{ab}=\widetilde{\nu}_a(x)\widetilde{\nu}_b(x)$ so that $\nabla_x$ is  $\Sp(2g,\Z)$ invariant.
\hfill \qed 

\medskip

Note that a $\Sp(2g,\Z)$ modular derivative can also be naturally defined generalizing the Serre derivative for $\SL(2,\Z)$ modular forms as follows.
Let $F_k(\Omega)$ denote a meromorphic $\Sp(2g,\Z)$ Siegel modular form of weight $k$ i.e. 
$F_k(\Omega)$ is meromorphic on $\HH_g$ where for all  $\gamma=\left( \begin{smallmatrix} A & B \\ C & D \end{smallmatrix} \right) \in\Sp(2g,\Z)$ 
\begin{align*}
F_k\left(\widetilde{\Omega}\right)=\det (C\Omega+D)^kF_k(\Omega).
\end{align*}
For projective connection $s(x)$  we define the projective differential $2$--form
\begin{align}
G_k(x,\Omega):=
\left(\nabla_x+\frac{k}{6}s(x)\right)F_k(\Omega).
\label{eq:modder}
\end{align}
From  \eqref{eq:projcmod} and Proposition~\ref{prop:nablaxmod} we have 
\begin{lemma}\label{lem:modder}
$G_k(x,\Omega)$  transforms under $\Sp(2g,\Z)$ like a weight $k$ Siegel modular form. 
\end{lemma}
\noindent 
This result  extends to a Siegel modular form  $F_k$ for a subgroup of $\Sp(2g,\Z)$ with a multiplier system.  
\medskip

\eqref{eq:modder} is a genus two version of the Serre modular derivative $g_{k+2}(q)=(q\partial_q +kE_2(q))f_k(q)$ for an $\SL(2,\Z)$ modular form $f_k(q)$ of  weight $k$ for which $g_{k+2}(q)$  is a modular form of weight $k+2$.
Equivalently, employing the standard $z$--coordinate on the torus,  $g_{k+2}(q)dz^2 $ transforms like a modular form of weight $k$ under $\SL(2,\Z)$. 

\subsection{Genus two surfaces formed from sewn tori}
We now review a general method due to Yamada \cite{Y} and discussed in detail in \cite{MT2} for calculating $\omega(x,y)$, $\nu_i(x)$ and $\Omega_{ij}$ for $i,j=1,2$ on the genus two Riemann surface formed by sewing together two tori $\mathcal{S}_{a} $ for $a=1,2$. We  sometimes refer to $\mathcal{S} _{1}$ and $\mathcal{S}_{2}$ as the left and right torus respectively.

Consider an oriented torus $\mathcal{S}_{a}=\C/\Lambda _{a}$ with lattice $\Lambda _{a}=2\pi i(\Z\tau _{a}\oplus \Z)$ for $\tau _{a}\in \HH_{1}$. For local coordinate $z_{a}\in \C/\Lambda _{a}$  the closed disc $\left\vert z_{a}\right\vert \leq r_{a}$  is contained in $\mathcal{S}_{a}$ provided $r_{a}<\frac{1}{2} D(q_{a})$ where
\begin{align*}
D(q_{a})=\min_{\lambda \in \Lambda _{a},\lambda \neq 0}|\lambda |,
\end{align*}
is the minimal lattice distance. We introduce a complex sewing   parameter $\epsilon$ where $|\epsilon |\leq r_{1}r_{2}<\frac{1}{4}D(q_{1})D(q_{2})$ and excise the disc $\{z_{a},\left\vert z_{a}\right\vert \leq |\epsilon |/r_{\bar{a}}\}$ centered at $z_{a}=0$ to form a punctured torus 
\begin{align}\label{Sahat}
\widehat{\mathcal{S}}_{a}=
\mathcal{S}_{a}\backslash\left \{z_{a},\left\vert z_{a}\right\vert \leq |\epsilon |/r_{\bar{a}}\right\},
\end{align}
where we here (and below) we use the convention 
\begin{align}\label{abar}
\overline{1}=2,\quad \overline{2}=1.
\end{align}

\noindent
Defining the annulus $
\mathcal{A}_{a}=\{z_{a},|\epsilon |/r_{\bar{a}}\leq \left\vert
z_{a}\right\vert \leq r_{a}\}$
we identify $\mathcal{A}_{1}$ with $\mathcal{A}_{2}$ via the sewing relation 
\begin{align}\label{pinch}
z_{1}z_{2}=\epsilon.  
\end{align}
The genus two Riemann surface $\mathcal{S}^{(2)}$ is parameterized by the sewing  domain\footnote{The sewing domain $\mathcal{D}_{\mathrm{sew}}$ is notated by $\mathcal{D}_{\epsilon}$ in \cite{MT2,MT3}} 
\begin{align}\label{Deps}
\mathcal{D}_{\mathrm{sew}}=\left\{(\tau _{1},\tau _{2},\epsilon )\in \HH_{1}
\mathbb{\times H}_{1}\mathbb{\times C}: |\epsilon |<\frac{1}{4}
D(q_{1})D(q_{2})\right\}.
\end{align}
We next introduce the infinite dimensional matrices 
\begin{align}
A_{a}= A_{a}(k,l,\tau _{a},\epsilon ) 
= \frac{ (-1)^{k+1}\epsilon ^{(k+l)/2}}{\sqrt{kl}} \frac{(k+l-1)!}{(k-1)!(l-1)!}E_{k+l}(\tau_a).
\label{eq:Aadef}
\end{align}
$A_1,A_2$ play an important role for the Heisenberg VOA on a genus two Riemann surface.  
Let $\mathbbm{1}$ denote  the infinite identity matrix and define
\begin{align*}
(\mathbbm{1}-A_{1}A_{2})^{-1}&=\sum_{n\geq 0}(A_{1}A_{2})^{n},\\
\log \det (\mathbbm{1}-A_{1}A_{2}) &=\mathrm{Tr}\log (\mathbbm{1}-A_{1}A_{2})   =-\sum_{n\geq 1}\frac{1}{n}\mathrm{Tr}((A_{1}A_{2})^{n}). 
\end{align*}
\begin{theorem}[\cite{MT2}] \label{Theorem_A1A2}
For all $(\tau _{1},\tau _{2},\epsilon )\in \mathcal{D}_{\mathrm{sew}}$ the matrix $(\mathbbm{1}-A_{1}A_{2})^{-1}$
is convergent and 
$\det (\mathbbm{1}-A_{1}A_{2})$ is non-vanishing and holomorphic.
\end{theorem}

\noindent
The bidifferential form $\omega(x,y)$, the holomorphic  1--differentials $\nu_i(x)$ and the period matrix $\Omega_{ij}$ are given in terms of the matrices $A_a$ and holomorphic one differentials on the punctured torus $\widehat{\mathcal{S}}_{a}$ defined by
\begin{align}
a(x;k)=\sqrt{k}\epsilon ^{k/2}P_{k+1}(x,\tau _{a})dx,
\label{eq:adef}
\end{align}
for $x\in \widehat{\mathcal{S}}_{a}$.
Letting $a(x)$, $a^T(x)$ denote the infinite row and column vector indexed by $k\ge 1$ we find
\begin{theorem}[\cite{MT2}] \label{Theorem_om_eps}
The genus two bidifferential  form $\omega(x,y)$ and the holomorphic  1--differentials $\nu_a(x)$for $a=1,2$ are given by
\begin{align*}
\omega(x,y)=&
\begin{cases}
P_2(x-y, \tau_a)dxdy
+ a(x)A_{\bar{a}}(\mathbbm{1}-A_a A_{\bar{a}})^{-1}a^{T}(y),
& x, y\in \widehat{\mathcal{S}}_{a}, \\ 
-a(x)(\mathbbm{1}-A_{\bar{a}}A_a)^{-1}a^{T}(y), 
& x\in \widehat{\mathcal{S}}%
_{a},\,y\in \widehat{\mathcal{S}}_{\bar{a}},
\end{cases}\\
\nu_a(x)=&
\begin{cases}
dx+\epsilon^{1/2}(a(x)A_{\bar{a}}(\mathbbm{1}-A_a A_{\bar{a}})^{-1})(1), 
& x\in \widehat{\mathcal{S}}_{a}, 
\\ 
-\epsilon^{1/2}(a(x)(\mathbbm{1}-A_a A_{\bar{a}})^{-1})(1), 
& x\in \widehat{\mathcal{S}}_{\bar{a}},
\end{cases}
\end{align*}
where $(1)$ refers to the $(1)$-entry of a vector.
\end{theorem} 
\begin{theorem}[\cite{MT2}]\label{Theorem_period_eps}
The sewing formalism
determines a holomorphic map 
\begin{eqnarray*}
F^{\Omega }:\mathcal{D}_{\mathrm{sew}} &\rightarrow &\HH_{2},  \notag
\\
(\tau _{1},\tau _{2},\epsilon ) &\mapsto &\Omega (\tau _{1},\tau
_{2},\epsilon ),  \label{Fepsmap}
\end{eqnarray*}%
where $\Omega =\Omega (\tau _{1},\tau _{2},\epsilon )$ is given by 
\begin{eqnarray*}
2\pi i\Omega _{11} &=&2\pi i\tau _{1}+\epsilon
(A_{2}(\mathbbm{1}-A_{1}A_{2})^{-1})(1,1),  \label{Om11eps} \\
2\pi i\Omega _{22} &=&2\pi i\tau _{2}+\epsilon
(A_{1}(\mathbbm{1}-A_{2}A_{1})^{-1})(1,1),  \label{Om22eps} \\
2\pi i\Omega _{12} &=&-\epsilon (\mathbbm{1}-A_{1}A_{2})^{-1}(1,1),
\end{eqnarray*}
where $(1,1)$ refers to the $(1,1)$--entry of a matrix. 
$F^{\Omega }$ is equivariant with respect to the action of $\Gamma\simeq (SL(2,\Z)\times SL(2,\Z))\rtimes \Z_{2}$ which preserves $\mathcal{D}_{\mathrm{sew}}$.
\end{theorem}
\noindent
Later we will  show below in Theorems~\ref{thm:Dxnabla1} and \ref{thm:Dxnabla2}  that $F^{\Omega }$ is injective but not surjective.


\section{Vertex Operator Algebras on Genus One and Two Riemann Surfaces}
\label{sec: VOA12review}

\subsection{Vertex operator algebras}
We review  aspects of vertex operator algebras e.g. \cite{B,FLM,FHL,Ka,LL,MN,MT5}. 
A Vertex Operator Algebra (VOA) is a quadruple $(V,Y,\mathbf{1},\omega)$ consisting of a $\Z$-graded complex vector space $V = \bigoplus_{n\in\Z}\,V_{(n)}$ where $\dim V_{(n)}<\infty$ for each $n\in \Z$, a linear map $Y:V\rightarrow\End(V)[[z,z^{-1}]]$ for a formal parameter $z$ and pair of distinguished vectors: the vacuum $\mathbf{1}\in V_{(0)}$ and the conformal vector $\omega\in V_{(2)}$. For each $v\in V$, the image under the map $Y$ is the vertex operator
\begin{align}
Y(v,z) = \sum_{n\in\Z}v(n)z^{-n-1},
\end{align}
with \emph{modes} $v(n)\in\End(V)$, where $Y(v,z)\mathbf{1} = v+O(z)$. Vertex operators satisfy \emph{locality} i.e. for all $u,v\in V$ there exists an integer $k\geq 0$ such that
\begin{align}
(z_1 - z_2)^k \left[ Y(u,z_1), Y(v,z_2) \right] = 0.
\end{align}
The vertex operator of the conformal vector $\omega$ is 
\begin{align*}
&Y(\omega,z) = \sum_{n\in\Z}L(n)z^{-n-2},
\end{align*}
where the modes $L(n)$ satisfy the Virasoro algebra with \emph{central charge} $c$
\begin{align}
&[L(m),L(n)] = (m-n)L(m+n) + \frac{m^3-m}{12}\delta_{m,-n}c\,\mathrm{Id}_V.
\label{eq:Ln}
\end{align}
We define the homogeneous space of weight $k$ to be 
\begin{align*}
V_{(k)} = \{v\in V | L(0)v = kv\},
\end{align*}
and we write $\wt(v)=k$ for $v\in V_{(k)}$. 
Finally, we have a translation condition
\begin{align}
&Y(L(-1)u,z) = \partial_z Y(u,z).
\label{eq:YT}
\end{align}

\bigskip
For a given VOA $V$, we define the \emph{adjoint vertex operator} (with respect to $A$)  by
\begin{align}
Y^\dagger (u,z) =\sum_{n\in\Z}u^\dagger(n) z^{-n-1}=  Y\left(\exp\left({\frac{z}{A}L(1)}\right)\left(- \frac{A}{z^2}\right)^{L(0)}u, \frac{A}{z}\right),
\label{eq:Ydag}
\end{align} 
associated with the formal M\"obius map $z\mapsto A/z$ \cite{FHL}.
For $u$ quasiprimary (i.e. $L(1)v=0$) of weight $\wt(u)$ then
\begin{align*}
u^\dagger(n) = (-1)^{\wt(u)} A^{n+1-\wt(u)} u(2\wt(u)-n-2).
\end{align*}
A bilinear form $\langle\, ,\, \rangle : V \times V\rightarrow \C$ is called \emph{invariant} if \cite{FHL,Li}
\begin{align}
\langle Y(u,z)a,b \rangle = \langle a,Y^\dagger(u,z)b \rangle\quad \mbox{ for all } a,b,u\in V.
\label{eq:bilform}
\end{align}
The adjoint vertex operator and $\langle\, ,\, \rangle$ depend on $A$. In particular
\begin{align}
\langle a,b\rangle\vert_{A=1}= A^{wt(a)}\langle a,b\rangle.
\label{eq:LiZA}
\end{align}
$\langle\, ,\, \rangle$  is necessarily symmetric \cite{FHL}.  In terms of modes, we have
\begin{align}
\langle u(n)a,b \rangle = \langle a,u^\dagger(n)b \rangle.
\label{eq:undag}
\end{align} 
Choosing $u=\omega$ and $n=1$ implies $\langle L(0)a,b \rangle = \langle a,L(0)b \rangle$. 
Thus $\langle a,b \rangle=0$ when $\wt(a) \neq \wt(b)$. 

A VOA is of \emph{strong CFT--type} if $V_{(0)} = \C\mathbf{1}$ and $V$ is simple and self-dual ($V$ is isomorphic to the dual module $V^\prime$ as a $V$-module). 
\cite{Li} guarantees that $V$ of strong CFT--type has a unique invariant non-degenerate bilinear form up to normalization. These results motivate the following definition
\begin{mydefn}\label{def:LiZ}
The \emph{Li-Zamolodchikov (Li-Z) metric}  on $V$ of strong CFT--type is the unique invariant bilinear form $\langle\,,\, \rangle $ normalized by $\langle \mathbf{1}, \mathbf{1} \rangle = 1$.
\end{mydefn}

\subsection{VOAs on a genus one Riemann surface}
Given a VOA $(V,Y,\mathbf{1},\omega)$, we can find an isomorphic VOA $(V,Y[,],\mathbf{1},\omt)$ introduced by Zhu \cite{Z}, called the \emph{square bracket} VOA. Both VOAs have the same underlying vector space $V$, vacuum vector $\mathbf{1}$ and central charge. The operators $Y[\,,\,]$ are defined by the coordinate change
\begin{align*}
Y[v,z] = \sum_{n\in \Z} v[n] z^{-n-1} = Y\left(q_z^{L(0)}v, q_z-1\right).
\end{align*} 
The new conformal vector is 
$
\omt = \omega - \frac{c}{24}\mathbf{1}
$,
with vertex operator
$
Y[\omt,z] = \sum_{n\in \Z} L[n] z^{-n-2}
$.
$L[0]$ provides an alternative $\Z$--grading on $V$ and we write  $\wt[v]=k$ if $L[0]v = kv$ where $\wt[v]=\wt(v)$ for primary $v$ ($L(n)v=0$ for all $ n>0$). We can similarly define a square bracket Li-Z metric $\langle\,, \,\rangle_{\mathrm{sq}}$.
The subscript sq will  be omitted where there is no ambiguity. 
\medskip

The \emph{genus one partition function} for $V$ is defined by the formal trace function
\begin{align*}
Z^{(1)}_V(\tau) = \Tr_V\left( q^{L(0)-c/24} \right),
\end{align*}
and the \emph{genus one $n$-point correlation function} is the formal expression 
\begin{align*}
Z^{(1)}_V(v_1,z_1;\ldots;v_n,z_n;\tau) = \Tr_V\left( Y(q_{z_1}^{L(0)}v_1,q_{z_1})\ldots Y(q_{z_n}^{L(0)}v_n,q_{z_n}) q^{L(0)-c/24} \right). 
\end{align*}
In particular, the genus one $1$-point function for $v\in V$ is
\begin{align*}
Z^{(1)}_V(v;\tau) &=  \Tr_V\left( o(v) q^{L(0)-c/24} \right),
\end{align*}
where, for $v$ homogeneous, $o(v) := v(\wt(v) - 1):V_{(m)}\rightarrow V_{(m)}$.
Every $n$-point function is expressible in terms of $1$-point functions \cite{MT1}
\begin{align}
&Z^{(1)}_V(v_1,z_1;\ldots;v_n,z_n;\tau)
\notag
\\
&= Z^{(1)}_V(Y[v_1,z_1]\ldots Y[v_n,z_n]\mathbf{1};\tau) 
\label{eq:Znpt1pt}
\\
&= Z^{(1)}_V(Y[v_1,z_{1}-z_{n}]\ldots Y[v_n,z_{n-1}-z_{n}]v_n;\tau).
\label{eq:Znpttrans}
\end{align}
We will make repeated use of Zhu recursion  which recursively relates formal $n$-point correlation functions to $(n-1)$-point functions \cite{Z} 
\begin{theorem}\label{zhu}[Zhu Recursion]
Genus one $n$-point correlation functions obey
\begin{align}
&Z^{(1)}_V(v_1,z_1;\ldots;v_n,z_n;\tau)
\notag
\\
&= \Tr_V \left( o(v_1) Y(q_{z_2}^{L(0)}v_2,q_{z_2})\ldots Y(q_{z_n}^{L(0)}v_n,q_{z_n}) q^{L(0)-c/24} \right)
\notag
 \\
&+ \sum_{k=2}^n \sum_{j\geq0} P_{1+j}(z_1-z_k,\tau) Z^{(1)}_V(v_2,z_2;\ldots;v_1[j]v_k,z_k;\ldots;v_n,z_n;\tau),
\label{eq:Zhurec}
\end{align}
for Weierstrass $P_k(z,\tau)$ of Definition~\ref{def:Pfunctions}.
\end{theorem}
\noindent
One of core ideas in Zhu theory is that the universal coefficients in the formal recursion formula \eqref{eq:Zhurec} are analytic Weierstrass functions. This ultimately implies convergence and modular properties for genus one partition and $n$--point functions for suitable VOAs \cite{Z}. 
The main aim of this paper is to find a genus two version of Theorem~\ref{zhu}.

Theorem~\ref{zhu} has many important applications e.g.   
$o(\omt) = L(0) - c/24$ implies
\begin{align}
Z^{(1)}_V(\omt;\tau) = q\partial_q Z^{(1)}_V(\tau).
\label{eq:delqZ}
\end{align}
For $n$ primary vectors $v_1,\ldots,v_n \in V$ one finds the genus one \emph{Ward Identity}
\begin{align}
&Z^{(1)}_V(\omt,x;v_1,x_1;\ldots;v_n,x_n;\tau)
\notag
\\
&= \left(q\partial_q + \sum_{k=1}^n \left(P_1(x-x_k,\tau)\partial_{x_k} + \wt[v_k]P_2(x-x_k,\tau) \right)\right) Z^{(1)}_V (v_1,z_1;\ldots;v_n,z_n;\tau).
\label{eq:Z1Ward}
\end{align}
We  also have that the Virasoro $n$-point function
\begin{align}
&Z^{(1)}_V(\omt,x_1;\ldots;\omt,x_n;\tau)
\notag
\\
&= \left(q\partial_q + \sum_{k=2}^n \left(P_1(x_1-x_k,\tau)\partial_{x_k} + 2P_2(x_1-x_k,\tau) \right)\right) Z^{(1)}_V(\omt,x_2;\ldots;\omt,x_n;\tau)
\notag
\\
&+\frac{c}{2}\sum_{k=2}^n P_4(x_1-x_k,\tau) Z^{(1)}_V(\omt,x_2;\ldots;\widehat{\omt,x_k};\ldots;\omt,x_n;\tau),
\label{eq:Z1Ward2}
\end{align}
where the caret indicates that the insertion of $\omt$ at $x_k$ is omitted. Theorem \ref{zhu} gives identities between formal series, while the differential equations from the Ward Identities allow us to prove convergence in many cases \cite{Z}.   
\begin{remark}
\label{rem:Mmodule}
The above definitions can be naturally extended to define $Z^{(1)}_{M}(\ldots)$ for a  graded $V$--module $M$ where  the trace is taken over $M$. 
\end{remark}
\begin{remark}
\label{rem:g1MDE} Modular differential equations for the genus one partition function can be obtained from \eqref{eq:Z1Ward} and \eqref{eq:Z1Ward} for any Virasoro singular vectors.
\end{remark}

\subsection{VOAs on a genus two Riemann surface}
We next review the definition of the formally defined genus two partition function and   $n$-point correlation functions for a VOA \cite{MT3}. 
We assume that $V$ is of strong CFT--type and hence a non-degenerate Li-Z metric exists. 
For a $V$--basis $\{u^{(a)}\}$ we define the dual basis $\{\overline{u}^{(a)}\}$ with respect to the Li-Z metric where
\begin{align*}
\langle u^{(a)}, \overline{u}^{(b)} \rangle_{\mathrm{sq}} = \delta_{ab}.
\end{align*}
The \emph{genus two partition function} for $V$ is defined by
\begin{align}
Z^{(2)}_V(\tau_1,\tau_2,\epsilon) = \sum_{u\in V} Z^{(1)}_V(u;\tau_1) Z^{(1)}_V(\overline{u};\tau_2),
\label{eq:Zdef}
\end{align}
where the formal sum is taken over the $V$--basis and $\overline{u}$ is the dual of $u$ with respect to $\langle,\rangle_{\mathrm{sq}}~$ with $A=\epsilon$ in \eqref{eq:Ydag}, i.e. we define the adjoint by
\begin{align*}
Y_\epsilon^\dagger[v,z]=Y\left[\exp\left({\frac{z}{\epsilon}L[1]}\right)\left(- \frac{\epsilon}{z^2}\right)^{L[0]}v, \frac{\epsilon}{z}\right].
\end{align*}
The \emph{genus two} $n$\emph{-point correlation function} for $a_1,\ldots,a_L$ 
and $b_1,\ldots,b_R$ formally  inserted at $x_1,\ldots,x_L\in\widehat{\mathcal{S}}_1$ and $y_1,\ldots,y_R\in\widehat{\mathcal{S}}_2$, respectively, by
\begin{align}
&Z^{(2)}_V(a_1,x_1;\ldots;a_L,x_L|b_1,y_1;\ldots;b_R,y_R;\tau_1,\tau_2,\epsilon) 
\notag
\\
&= \sum_{u\in V} Z^{(1)}_V(Y[a_1,x_1]\ldots Y[a_L,x_L]u;\tau_1) Z^{(1)}_V(Y[b_R,y_R]\ldots Y[b_1,y_1]\overline{u};\tau_2),
\label{eq:Znptdef}
\end{align}
where the sum as in \eqref{eq:Zdef}.
\begin{remark}
\label{rem:basis}
\eqref{eq:Zdef} and \eqref{eq:Znptdef} are independent of the choice of $V$--basis.
\end{remark}
\begin{remark}
\label{rem:M1M2}
As with Remark~\ref{rem:Mmodule}, the above definitions can be  extended to define $Z^{(2)}_{M_1M_2}(\ldots)$ for a  pair of $V$--modules $M_1,M_2$, where the left and right--hand genus one contributions in \eqref{eq:Zdef} or \eqref{eq:Znptdef} are replaced by trace functions over $M_1$ and  $M_2$, respectively.
\end{remark}

\begin{remark}
\label{rem:ZMTdef}
\eqref{eq:Zdef} is equivalent to the original definition of  \cite{MT3} where the the $\epsilon$ dependence is made explicit:
\begin{align*}
Z^{(2)}_V(\tau_1,\tau_2,\epsilon) = \sum_{r\geq0} \epsilon^r \sum_{u\in V_{[r]}} Z^{(1)}_V(u;\tau_1) Z^{(1)}_V(\overline{u};\tau_2),
\end{align*}
and the internal sum is taken over any $V_{[n]}$--basis and $\overline{u}$ is the dual of $u$ with respect to the Li-Z metric and adjoint operators defined by the mapping $z\mapsto 1/z$. 
Definition  \eqref{eq:Zdef} has the benefit of streamlining later analysis. Similar remarks apply to \eqref{eq:Znptdef}.
\end{remark}



\section{Zhu Reduction for Genus Two $n$-Point Correlation Functions}
\label{sec:Zhured}
\subsection{Genus two $n$-point correlation functions}
In this section we derive a formal genus two Zhu reduction expression for all  $n$-point correlation functions. Let  $v\in V$ be inserted at $x\in\widehat{\mathcal{S}}_1$, $a_1,\ldots,a_L\in V$ be inserted at $x_1,\ldots,x_L\in\widehat{\mathcal{S}}_1$ and 
$b_1,\ldots,b_R\in V$ be inserted at $y_1,\ldots,y_R\in\widehat{\mathcal{S}}_2$.
We consider  the corresponding genus two $n$-point function
\begin{align}
&Z^{(2)}_V (v,x;\boldsymbol{a_l,x_l}|\boldsymbol{b_r,y_r};\tau_1, \tau_2, \epsilon) \notag \\
&= \sum_{u\in V} 
Z^{(1)}_V(Y[v,x]\boldsymbol{Y[a_l,x_l]}u;\tau_1)
Z^{(1)}_V(\boldsymbol{Y[b_r,y_r]}\overline{u};\tau_2),
\label{Znvleft1}
\end{align}
with the following notational abbreviations:
\begin{align*}
&\boldsymbol{a_l,x_l}\equiv a_1,x_1;\ldots;a_lx_L,\qquad 
\boldsymbol{Y[a_l,x_l]}\equiv Y[a_1,x_1]\ldots Y[a_L,x_L],
\\
& 
\boldsymbol{b_r,y_r}\equiv b_1,y_1;\ldots;b_R,y_R,\qquad
\boldsymbol{Y[b_r,y_r]}\equiv Y[b_1,y_1]\ldots Y[b_R,y_R].
\end{align*}
There is a similar expression for $x\in\widehat{\mathcal{S}}_2$ with the $Y[v,x]$ vertex operator  inserted on the right hand side of \eqref{Znvleft1}.
Zhu recursion (Theorem~\ref{zhu}) implies
\begin{align}
&Z^{(1)}_V(Y[v,x]\boldsymbol{Y[a_l,x_l]}u;\tau_1) \notag \\
&= \Tr_V \left(o(v)\boldsymbol{Y(q_{x_l}^{L(0)}a_l,q_{x_l})} Y(q_{0}^{L(0)}u,q_{0})q_1^{L(0) - c/24} \right) \notag \\
&+ \sum_{l=1}^L \sum_{j\geq 0} P_{1+j}(x-x_l,\tau_1)Z^{(1)}_V (\ldots ;v[j]a_l,x_l;\ldots;\tau_1) \notag \\
&+ \sum_{m\geq 0} P_{1+m}(x,\tau_1) Z^{(1)}_V(\boldsymbol{Y[a_l,x_l]}v[m]u;\tau_1),
\label{Znvleft2}
\end{align}
where $\boldsymbol{Y(q_{x_l}^{L(0)}a_l,q_{x_l})}\equiv Y(q_{x_1}^{L(0)}a_1,q_{x_1})\ldots Y(q_{x_L}^{L(0)}a_L,q_{x_L})$ and $q_a = e^{\tpi \tau_a}$.

To streamline notation, we make a number of definitions. 
We will often suppress the explicit dependence on $v,\boldsymbol{a_l,x_l},\boldsymbol{b_r,y_r}$ and $\tau_1, \tau_2, \epsilon$ when there is no ambiguity.
Let $O_a$   for $a\in \{1,2\}$  be defined  by
\begin{align}
&O_1 = O_1(v;\boldsymbol{a_l,x_l}|\boldsymbol{b_r,y_r};\tau_1, \tau_2, \epsilon) \notag \\
&= \sum_{u\in V} \Tr_V \left( o(v)\boldsymbol{Y(q_{x_l}^{L(0)}a_l,q_{x_l})} Y(q_{0}^{L(0)}u,q_{0})q_1^{L(0) - c/24} \right) Z^{(1)}_V(\boldsymbol{Y[b_r,y_r]}\overline{u};\tau_2),
\notag
 \\
&O_2 = O_2(v;\boldsymbol{a_l,x_l}|\boldsymbol{b_r,y_r};\tau_1, \tau_2, \epsilon) \notag \\
&= \sum_{u\in V} Z^{(1)}_V(\boldsymbol{Y[a_l,x_l]}u;\tau_1) \Tr_V \left( o(v)\boldsymbol{Y(q_{y_r}^{L(0)}b_r,q_{y_r})} Y(q_{0}^{L(0)}\overline{u},q_{0})q_2^{L(0) - c/24} \right),
\label{eq:Oadef}
\end{align}
where   $q_{0}=1$ using the translation property \eqref{eq:Znpttrans}. 
 
We define a number of infinite matrices and row and column vectors indexed by $m,n\ge 1$ as follows. Let $\Lambda_a $  for $a\in \{1,2\}$ be  the matrix with components
\begin{align}
\Lambda_a (m,n) = \Lambda_a (m,n;\tau_a,\epsilon) = 
\epsilon^{(m+n)/2}(-1)^{n+1} \binom{m+n-1}{n}E_{m+n}(\tau_a).
\label{eq:Lambda}
\end{align}
Note that  
\begin{align}
\Lambda_a =SA_{a}S^{-1}, 
\label{eq:LamA}
\end{align}
for $A_a$ of \eqref{eq:Aadef} for $S$ a diagonal matrix with components
\begin{align}
S(m,n)=\sqrt{m}\delta_{mn}.
\label{eq:Sdef}
\end{align}
Let $\mathbb{R}(x) $ for $x\in \widehat{\mathcal{S}}_a$ be the row vector with components
\begin{align}
\mathbb{R}(x;m) = \epsilon^{\frac{m}{2}} P_{m+1} (x, \tau_a) .
\label{eq:Rdef}
\end{align}
Let $\mathbb{X}_a $   for $a\in \{1,2\}$  be the column vector with components
\begin{align}
&\mathbb{X}_1(m) = 
\mathbb{X}_1(m;v;\boldsymbol{a_l,x_l}|\boldsymbol{b_r,y_r};\tau_1, \tau_2, \epsilon) 
\notag
\\
&= \epsilon^{-m/2}\sum_{u\in V} Z^{(1)}_V(\boldsymbol{Y[a_l,x_l]}v[m]u;\tau_1) Z^{(1)}_V(\boldsymbol{Y[b_r,y_r]}\overline{u};\tau_2), 
\notag
\\
&\mathbb{X}_2(m) = 
\mathbb{X}_2(m;v;\boldsymbol{a_l,x_l}|\boldsymbol{b_r,y_r};\tau_1, \tau_2, \epsilon) 
\notag
\\
&= \epsilon^{-m/2}\sum_{u\in V} Z^{(1)}_V(\boldsymbol{Y[a_l,x_l]}u;\tau_1) Z^{(1)}_V(\boldsymbol{Y[b_r,y_r]}v[m]\overline{u};\tau_2).
\label{eq:Xadef}
\end{align}
Lastly, define genus two contraction terms for $j\ge 0$ given by
\begin{align}
Z^{(2)}_V(\ldots ;v[j]a_l,x_l; \ldots) = \sum_{u\in V}Z^{(1)}_V(\ldots v[j]a_l,x_l \ldots;\tau_1) Z^{(1)}_V(\boldsymbol{Y[b_r,y_r]}\overline{u};\tau_2), 
\notag\\
Z^{(2)}_V(\ldots ;v[j]b_r,y_r; \ldots) = \sum_{u\in V}Z^{(1)}_V(\boldsymbol{Y[a_l,x_l]}u;\tau_1)Z^{(1)}_V(\ldots v[j]b_r,y_r \ldots;\tau_2).
\label{eq:Z2vj}
\end{align}
Thus applying genus 1 Zhu reduction \eqref{Znvleft2} to \eqref{Znvleft1}  using \eqref{eq:Oadef}--\eqref{eq:Z2vj}  we find
\begin{align}
&Z^{(2)}_V (v,x;\boldsymbol{a_l,x_l}|\boldsymbol{b_r,y_r}) 
= O_1 
+ \mathbb{R}(x)\mathbb{X}_1 
\notag \\
&+ \sum_{l=1}^L \sum_{j\geq 0} P_{1+j}(x-x_l,\tau_1)Z^{(2)}_V (\ldots ;v[j]a_l,x_l;\ldots) 
\notag \\
&+ P_1(x,\tau_1) \sum_{u\in V}Z^{(1)}_V(\boldsymbol{Y[a_l,x_l]}v[0]u;\tau_1) Z^{(1)}_V(\boldsymbol{Y[b_r,y_r]}\overline{u};\tau_2),
\label{Znvleft3}
\end{align}
and  similarly  for $v$ inserted on the right hand side of \eqref{Znvleft1} with $x\in\widehat{\mathcal{S}}_2$.

\subsection{A recursive identity for $\mathbb{X}_a$ }
We next  develop a recursive formula for $\mathbb{X}_a$ of \eqref{eq:Xadef}
which can  be formally solved to obtain a closed expression for $n$-point functions \eqref{Znvleft1} in terms of $n-1$ point functions with universal coefficients. 
Assume that $v$ is quasiprimary of weight $wt[v]=N$
(we consider quasiprimary descendants later).
Let $\langle \,,\, \rangle$
denote the square bracket Li-Z metric of \eqref{eq:Ydag} with $A = \epsilon$. 
Then using \eqref{eq:undag} we find
\begin{align*} 
v[m]u &= \sum_{w\in V} \left\langle \overline{w}, v[m]u \right\rangle w 
= \sum_{w\in V} \langle  v^{\dagger}[m]\overline{w}, u \rangle w\\
& = (-1)^N \epsilon^{m-K/2} \sum_{w\in V}\left\langle  v[K-m]\overline{w}, u \right\rangle w ,
\end{align*}
where the $w$ sum is taken over any $V$--basis and where 
\begin{align}
K=2(N-1).
\label{eq:Kdef}
\end{align}  
Since $\sum_{u\in V} \left\langle  v[K-m]\overline{w}, u \right\rangle \overline{u}=v[K-m]\overline{w}$ we thus find
\begin{align}
&\mathbb{X}_1(m) 
= (-1)^N \epsilon^{(m-K)/2}\sum_{w\in V} Z^{(1)}_V(\boldsymbol{Y[a_l,x_l]}w;\tau_1) Z^{(1)}_V(\boldsymbol{Y[b_r,y_r]}v[K-m]\overline{w};\tau_2)
\label{nX1_1}.
\end{align}
Provided $m\ge K+1$ then genus one Zhu recursion applies to the right hand side of \eqref{nX1_1} leading to a recursive identity for $\mathbb{X}_1(m)$. 
However,  provided $K\ge 2$ (i.e. $N\ge 2$) then the first $K$ components of $\mathbb{X}_1$ are not subject to this recursive formula.  
Zhu recursion  implies that for $s \geq 1$ \cite{Z}
\begin{align*}
&Z^{(1)}_V(\boldsymbol{Y[b_r,y_r]}v[-s]\overline{w};\tau_2) \\
&=\delta_{1s} \Tr_V\left(o(v)\boldsymbol{Y(q_{y_r}^{L(0)}b_r,q_{y_r})}Y(q_{0}^{L(0)}\overline{w},q_{0})q_2^{L(0)-c/24}\right) \\ 
&+\sum_{j\geq 0} (-1)^{j+1} \binom{s+j-1}{j} E_{s+j}(\tau_2) Z^{(1)}_V(\boldsymbol{Y[b_r,y_r]}v[j]\overline{w};\tau_2)\\
&+\sum_{r=1}^R \sum_{j\geq0} (-1)^{s+1}\binom{s+j-1}{j} P_{s+j}(y_r,\tau_2) Z^{(1)}_V(\ldots Y[v[j]b_r,y_r] \ldots\overline{w};\tau_2) .
\end{align*}
But Proposition~4.3.1 of \cite{Z} implies
\begin{align*}
Z^{(1)}_V(\boldsymbol{Y[b_r,y_r]}v[0]\overline{w};\tau_2) 
= -\sum_{r=1}^R Z^{(1)}_V(\ldots Y[v[0]b_r,y_r]\ldots \overline{w};\tau_2),
\end{align*} 
so that  using \eqref{eq:Lambda} 
\begin{align}
&\epsilon^{s/2}Z^{(1)}_V(\boldsymbol{Y[b_r,y_r]}v[-s]\overline{w};\tau_2) \notag\\
&=\delta_{1s} \epsilon^{1/2}\Tr_V\left(o(v)\boldsymbol{Y(q_{y_r}^{L(0)}b_r,q_{y_r})}Y(q_{0}^{L(0)}\overline{w},q_{0})q_2^{L(0)-c/24}\right) \notag\\ 
&+\sum_{j\geq 1} \Lambda_2(s,j) \epsilon^{-j/2}Z^{(1)}_V(\boldsymbol{Y[b_r,y_r]}v[j]\overline{w};\tau_2)\notag \\
&+\sum_{r=1}^R \sum_{j\geq0} (-1)^{j+1}
\mathbb{P}_{1+j}(y_r;s)
Z^{(1)}_V(\ldots Y[v[j]b_r] \ldots\overline{w};\tau_2),
\label{eq:ZVs}
\end{align}
where $\mathbb{P}_{1+j}(x)=\frac{(-1)^j}{j!}\mathbb{P}_{1}(x)$, for $x\in \widehat{\mathcal{S}}_a$ and $j\ge 0$, is the column vector 
with components
\begin{align}
\mathbb{P}_{1+j}(x;m)=\epsilon^{\frac{m}{2}}\binom{m+j-1}{j}
\left( P_{m+j}(x,\tau_a) - \delta_{j 0}E_m(\tau_a)\right). 
\label{eq:P1jdef}
\end{align}
Substituting \eqref{eq:ZVs} into \eqref{nX1_1} with $s=m-K$ we therefore find that for $m\ge K+1$ 
\begin{align}
\mathbb{X}_1(m) 
=& (-1)^N \epsilon^{1/2}\delta_{m-K,1}O_2
+  (-1)^N \left( \Lambda_2\mathbb{X}_2 \right)(m-K)
\notag \\
&+ (-1)^N \sum_{r=1}^R\sum_{j\geq 0} (-1)^{j+1} 
\mathbb{P}_{1+j}(y_r;m-K) Z^{(2)}_V(\ldots ;v[j]b_r,y_r;\ldots).
\label{nX1_2}
\end{align}
Thus for $m\ge K+1$,  we can recursively relate $\mathbb{X}_1(m)$ to the $m-K$ component  of an infinite vector involving $\mathbb{X}_2$. In order to describe this index translation by $K = 2N - 2\ge 0$ we define infinite matrices
$\Gamma $ and $ \Delta $  with components
\begin{align}
\Gamma(m,n) &= \delta_{m, -n+K}, \quad
\Delta(m,n) = \delta_{m, n+K}.
\label{eq:GamDelTh}
\end{align}
We also define the projection matrix
\begin{eqnarray}
\Pi =\Gamma^2=
\begin{bmatrix}
\mathbbm{1}_{K-1} & 0 \\
0 & \ddots 
\end{bmatrix}
,
\label{eq:PiK}
\end{eqnarray}
where $\mathbbm{1}_{K-1}$ denotes the $K-1$ dimensional identity matrix (with $\mathbbm{1}_{-1}=0$).

\begin{lemma}\label{lem:matrices}
The matrices $\Gamma,\Delta,\Pi$ obey the identities
\begin{align}
\Gamma=\Pi \Gamma=\Gamma\Pi ,\quad\Gamma \Delta=\Delta^T\Gamma = 0,\quad \Delta^T\Delta = \mathbbm{1}.
\label{eq:Matrices}
\end{align}
\end{lemma}
\noindent We define  column  vectors $\mathbb{O}_{a}$ with one non-zero component
\begin{align*}
\mathbb{O}_{a}(m) = \epsilon^{1/2}\delta_{1m}O_{a},\quad a=1,2.
\end{align*}
(suppressing dependence on $v,\boldsymbol{a_l,x_l},\boldsymbol{b_r,y_r}$ etc.)
and   column vectors $\mathbb{G}_{a}$ given by
\begin{align}
\mathbb{G}_{1}
&= \sum_{l=1}^L\sum_{j\geq0} (-1)^{j+1} 
\mathbb{P}_{1+j}(x_l)
 Z^{(2)}_V(\ldots;v[j]a_l,x_l;\ldots),
\notag
\\
\mathbb{G}_{2}
&= \sum_{r=1}^R \sum_{j\geq0} (-1)^{j+1}
\mathbb{P}_{1+j}(y_r)
 Z^{(2)}_V(\ldots;v[j]b_r,y_r;\ldots).
\label{eq:Ga}
\end{align}
\noindent Then we can rewrite \eqref{nX1_2} for $m\geq K+1$ as
\begin{align}
\mathbb{X}_1(m) = (-1)^N \Big( \Delta\left(\mathbb{O}_2  + \mathbb{G}_{2}+ \Lambda_{2}\mathbb{X}_2\right)\Big)(m). 
\label{eq:X1X2rel}
\end{align}
with a similar formula for $\mathbb{X}_2(m)$ in terms of $\mathbb{O}_1$, $\mathbb{X}_1$ and $\mathbb{G}_{1}$. 

For $N>1$  the remaining components of $\mathbb{X}_1(m)$ for $1\leq m\leq K $ are described as follows. We first note from \eqref{nX1_1} that  for $1\leq m\leq K-1$ 
\begin{align}
\mathbb{X}_a(m) = &(-1)^N\mathbb{X}_{\overline{a}}(K-m) 
=
(-1)^N \left( \Gamma \mathbb{X}_{\overline{a}}\right)(m),
\label{nX1GamX2}
\end{align} 
for $a=1,2$ (recalling the convention $\overline{1}=2$ and $\overline{2}=1$).
Define the projection on to the  first $K-1$ components of $\mathbb{X}_a$ by
\begin{align}
\mathbb{X}_a^\Pi = \Pi \mathbb{X}_a, 
\label{nXPi}
\end{align}
(where $\mathbb{X}_a^\Pi =0$ if $K=0$).
Using  Lemma~\ref{lem:matrices} we may  rewrite  \eqref{nX1GamX2} as 
\begin{align}
\mathbb{X}_a^{\Pi} &= (-1)^N \Gamma\mathbb{X}_{\overline a}= (-1)^N \Gamma\mathbb{X}_{\overline a}^{\Pi}.
\label{nX1GamX22}
\end{align}
\eqref{nX1_1} also implies
\begin{align}
\mathbb{X}_1(K) = (-1)^N \sum_{u\in V}Z^{(1)}_V(\boldsymbol{Y[a_l,x_l]}u;\tau_1)Z^{(1)}_V(\boldsymbol{Y[b_r,y_r]}v[0]\overline{u};\tau_2).
\label{eq:X1Kdef}
\end{align}
By Proposition~4.3.1 of \cite{Z} this can be re-expressed as
\begin{align}
\mathbb{X}_1(K) &= -(-1)^N \sum_{r=1}^R Z^{(2)}_V(\ldots ;v[0]b_r,y_r;\ldots),
\label{eq:X1K}
\end{align}
and similarly
\begin{align}
\mathbb{X}_2(K) &= -(-1)^N \sum_{l=1}^L Z^{(2)}_V(\ldots ;v[0]a_l,x_l;\ldots).
\label{eq:X2K}
\end{align} 
Introducing an infinite vector $\mathbb{X}_a^K = \left(\mathbb{X}_a^K(m)\right)$ with one non-zero component
\begin{align}
\mathbb{X}_a^K(m) = \delta_{mK}\mathbb{X}_a(K),
\label{eq:XaK}
\end{align}
we therefore find altogether that 
\begin{proposition}
Let $v$ be a quasi-primary vector with  $\wt[v]=N$. Then $\mathbb{X}_a$ obeys the  recursive identity 
\begin{align}
\mathbb{X}_a = &(-1)^N  \Gamma\mathbb{X}_{\overline{a}}^{\Pi}
+\mathbb{X}_a^K 
+ (-1)^N \Delta\Big( \mathbb{O}_{\overline{a}} +\mathbb{G}_{\overline{a}}
+\Lambda_{\overline{a}}\mathbb{X}_{\overline{a}}
\Big).
\label{nXarec}
\end{align}
\end{proposition}
We next describe how to formally solve the recursive identity \eqref{nXarec}.
Let 
\begin{align*}
\mathbb{X}^{\perp}_a =\Delta^T\mathbb{X}_a.
\end{align*}
Decompose $\mathbb{X}_a$ as
\begin{align}
\mathbb{X}_a =\mathbb{X}^{\Pi}_a + \mathbb{X}_a^K + \Delta\mathbb{X}^{\perp}_a.
\label{eq:Xdecomp}
\end{align}
Since $\Delta^T\Gamma = 0$ and $\Delta^T\Delta = \mathbbm{1}$  it follows from \eqref{nXarec} and \eqref{eq:Xdecomp} that 
\begin{align}
\mathbb{X}^{\perp}_a=& (-1)^N \left(\mathbb{O}_{\overline{a}} 
 +\mathbb{G}_{\overline{a}}
+ \Lambda_{\overline{a}}\mathbb{X}_{\overline{a}}\right)\notag
\\
=& (-1)^N \left(\mathbb{O}_{\overline{a}} 
+ \mathbb{G}_{\overline{a}}
+ \Lambda_{\overline{a}}\left(\mathbb{X}_{\overline{a}}^{\Pi} 
+ \mathbb{X}_{\overline{a}}^K \right)
+ \widetilde{\Lambda}_{\overline{a}}\mathbb{X}^{\perp}_{\overline{a}} 
\right) ,
\label{nXth_it}
\end{align}
where we define
\begin{equation}
\widetilde{\Lambda}_a = \Lambda_a \Delta. 
\label{eq:Lamdatilde}
\end{equation}
Iterating \eqref{nXth_it}  we find
\begin{align*}
\mathbb{X}^{\perp}_a=& (-1)^N \left(\mathbb{O}_{\overline{a}}
+\mathbb{G}_{\overline{a}} 
+ \Lambda_{\overline{a}}\left(\mathbb{X}_{\overline{a}}^{\Pi} + \mathbb{X}_{\overline{a}}^K \right) \right) \\
&+ \widetilde{\Lambda}_{\overline{a}} 
\left(\mathbb{O}_{a} +\mathbb{G}_{a}
+ \Lambda_{a}\left(\mathbb{X}_{a}^{\Pi} + \mathbb{X}_{a}^K \right) 
+\widetilde{\Lambda}_{a}\mathbb{X}^{\perp}_{a}
\right). 
\end{align*}
Thus we may formally solve for $\mathbb{X}^{\perp}_a$ 
in terms of the formal matrix inverse 
\begin{align}
\left( \mathbbm{1} - \widetilde{\Lambda}_{\overline{a}}\widetilde{\Lambda}_{a} \right)^{-1}
= \sum_{n\geq 0} \left( \widetilde{\Lambda}_{\overline{a}}\widetilde{\Lambda}_{a} \right)^n,
\label{eq:IminLamLaminv}
\end{align}
for $a=1,2$. We therefore find
\begin{proposition}
\label{prop:Xth}
Let $v$ be a quasi-primary vector with  $\wt[v]=N$. Then $\mathbb{X}_a =\mathbb{X}^{\Pi}_a + \mathbb{X}_a^K + \Delta\mathbb{X}^{\perp}_a$ where
\begin{align}
\mathbb{X}^{\perp}_a=
&~ (-1)^N \left(\mathbbm{1} 
-\widetilde{\Lambda}_{\overline{a}}\widetilde{\Lambda}_{a}\right)^{-1} \left(\mathbb{O}_{\overline{a}} + \mathbb{G}_{\overline{a}}
+ \Lambda_{\overline{a}}\left(\mathbb{X}_{\overline{a}}^{\Pi} +\mathbb{X}_{\overline{a}}^K \right)
\right) 
\notag \\
&~+ \left(\mathbbm{1} -\widetilde{\Lambda}_{\overline{a}}\widetilde{\Lambda}_{a}\right)^{-1} \widetilde{\Lambda}_{\overline{a}} 
\left(\mathbb{O}_{a} +  \mathbb{G}_{a}
+ \Lambda_{a}\left(\mathbb{X}_{a}^{\Pi} + \mathbb{X}_{a}^K \right) 
\right).
\label{nXth}
\end{align}
\end{proposition}

\medskip
\subsection{Genus two Zhu recursion}
We now return to the original genus two $n$-point function \eqref{Znvleft3}. 
Substituting $\mathbb{X}_1 $ from Proposition~\ref{prop:Xth} we obtain
\begin{align}
&Z^{(2)}_V (v,x;\boldsymbol{a_l,x_l}|\boldsymbol{b_r,y_r}) 
\notag \\ 
&= O_1 
\notag 
+ \mathbb{R}(x)\left(\mathbb{X}^\Pi_1 + \mathbb{X}^K_1\right) 
\notag \\ 
&+ (-1)^N. \sups{N}\mathbb{Q}(x) \left(\mathbb{O}_{2}  
+\mathbb{G}_{2}
+ \Lambda_{2}\left(\mathbb{X}_{2}^{\Pi} + \mathbb{X}_{2}^K\right)
\right)
\notag \\
&+ \sups{N}\mathbb{Q}(x) \widetilde{\Lambda}_{2} 
\left(\mathbb{O}_{1}  + \mathbb{G}_{1}
+  \Lambda_{1}\left(\mathbb{X}_{1}^{\Pi} + \mathbb{X}_{1}^K\right)
\right) 
\notag \\
&+ \sum_{l=1}^L \Big(\sum_{j\geq 0} P_{1+j}(x-x_l,\tau_1)Z^{(2)}_V (\ldots ;v[j]a_l,x_l;\ldots) 
\notag \\
&- P_1(x,\tau_1) Z^{(2)}_V(\ldots;v[0]a_l,x_l;\ldots)\Big),
\label{Zn}
\end{align}
where $\sups{N}\mathbb{Q}(x)$ is an infinite row vector defined  by
\begin{equation}
\sups{N}\mathbb{Q}(x) =\mathbb{R}(x) \Delta \left( \mathbbm{1} - \widetilde{\Lambda}_{\overline{a}}\widetilde{\Lambda}_a \right)^{-1},\quad \mbox{ for }x\in\widehat{\mathcal{S}}_a.
\label{eq:Qdef}
\end{equation}
The pre-superscript $N$ is introduced to emphasise the dependence of this expression on $N$ through  $\Delta$ (recalling also that $\widetilde{\Lambda}_a=\Lambda_a\Delta$). 
\medskip

We next identify various contributing terms to \eqref{Zn}. Using \eqref{nX1GamX22}, \eqref{eq:X1K}
and  \eqref{eq:X2K} we can describe the $O_a,\mathbb{X}_{a}^{\Pi}$ coefficients in terms of the following:
\begin{mydefn}\label{calFforms}
Let $\supsin{N}\mathcal{F}_{a}(x)$ for $N\ge 1$ and $a=1,2$  be given by
\begin{align}
&
\supsin{N}\mathcal{F}_{a}(x)=
\begin{cases}
1 + \epsilon^{1/2} \Big( \sups{N}\mathbb{Q}(x)\widetilde{\Lambda}_{\overline{a}} \Big) (1),
\quad \mbox{ for }x\in\widehat{\mathcal{S}}_a,\\
(-1)^N\epsilon^{1/2} \Big(\sups{N}\mathbb{Q}(x)  \Big)(1)
,\quad \mbox{ for }x\in\widehat{{\mathcal S}}_{\abar},
\end{cases}
\label{eq:calFadef}
\end{align}
and let
$\supsin{N}\mathcal{F}^\Pi(x)$, for $x \in {\mathcal S}_a$,  be an infinite row vector  
given by
\begin{align}
\supsin{N}\mathcal{F}^\Pi(x)=
\left( \mathbb{R}(x) + \sups{N}\mathbb{Q}(x) \left(\widetilde{\Lambda}_{\overline{a}}\Lambda_{a} + \Lambda_{\overline{a}}\Gamma \right) \right)\Pi.
\label{eq:calFPidef}
\end{align}
Note that $\supsin{1}\mathcal{F}^\Pi(x)=0$ and otherwise $\supsin{N}\mathcal{F}^\Pi(x;m)=0$ for $m\ge K=2N-2$.
\end{mydefn}
\noindent The $\mathbb{X}^K_a$ terms of \eqref{Zn}  contribute
\begin{align*}
&\left(\mathbb{R}(x)
+ \sups{N}\mathbb{Q}(x) \widetilde{\Lambda}_{2}\Lambda_{1}\right)\mathbb{X}_{1}^{K}=\\
&\pi_N(-1)^{N+1} \left(
\epsilon^{K/2}P_{K+1}(x)
+ \left(\sups{N}\mathbb{Q}(x) \widetilde{\Lambda}_{2}\Lambda_{1}\right)(K)
\right)
\sum_{r=1}^R Z^{(2)}_V(\ldots ;v[0]b_r,y_r;\ldots),\\
 &(-1)^N. \sups{N}\mathbb{Q}(x)\Lambda_2  \mathbb{X}_{2}^K
= -\pi_N\left(\sups{N}\mathbb{Q}(x) \Lambda_{2}\right)(K)
\sum_{l=1}^L Z^{(2)}_V(\ldots ;v[0]a_l,x_l;\ldots),
\end{align*} 
using \eqref{eq:X1Kdef} and \eqref{eq:XaK} and  where $\pi_N\equiv 1-\delta_{1N}$ for $N\ge 1$. 
There are further contributions to the multipliers of the contraction terms
$Z^{(2)}_V(\ldots;v[j]a_l,x_l;\ldots)$ and $Z^{(2)}_V(\ldots;v[j]b_r,y_r;\ldots)$ for $j\geq0$ arising from the $\mathbb{G}_{a}$ terms \eqref{eq:Ga} and the last summation term in \eqref{Zn}. These can be described as  follows:
\begin{mydefn}\label{calP21}
Define $\sups{N}\mathcal{P}_{1}(x,y)=\sups{N}\mathcal{P}_{1}(x,y;\tau_1,\tau_2,\epsilon)$ for $N\ge 1$ by
\begin{align*}
\sups{N}\mathcal{P}_{1}(x,y) 
=&
P_1(x-y,\tau_a)- P_1(x,\tau_a)  - \sups{N}\mathbb{Q}(x)\widetilde{\Lambda}_{\abar} \,\mathbb{P}_{1} (y)
-\pi_{N}
\left(\sups{N}\mathbb{Q}(x)\Lambda_{\abar}\right) (K),
\end{align*}
for   $x,y\in\widehat{\mathcal{S}}_a$ and  
\begin{align*}
\sups{N}\mathcal{P}_{1}(x,y) 
=&
(-1)^{N+1} \Big[\sups{N}\mathbb{Q}(x) \mathbb{P}_{1}  (y)
 +\pi_{N} \epsilon^{K/2}P_{K+1}(x)
+\pi_{N}  \left( \sups{N}\mathbb{Q}(x)\widetilde{\Lambda}_{\abar}\Lambda_a \right)(K)\Big], 
\end{align*}
for $ x\in\widehat{\mathcal{S}}_a,\;y\in\widehat{\mathcal{S}}_{\abar}$ where $\pi_N=1-\delta_{N1}$ and $K=2N-2$.
\end{mydefn}

\begin{mydefn}\label{calP21+j}
For $j> 0$ define $\sups{N}\mathcal{P}_{1+j}(x,y)
= \frac{1}{j!} \partial_y^j \left(\sups{N}\mathcal{P}_{1}(x,y) \right)$, i.e.
\begin{align}
\sups{N}\mathcal{P}_{1+j}(x,y)=
\begin{cases}
P_{1+j}(x-y)+ (-1)^{1+j}.\sups{N}\mathbb{Q}(x)\widetilde{\Lambda}_{\abar} \,\mathbb{P}_{1+j}(y),\; &\mbox{for }x,y\in\widehat{\mathcal{S}}_a,\\
(-1)^{N+1+j}.\sups{N}\mathbb{Q}(x) \mathbb{P}_{1+j}(y)
,\; &\mbox{for }x\in\widehat{\mathcal{S}}_a,\, y\in\widehat{\mathcal{S}}_{\abar}.
\end{cases}
\label{eq:PN21j}
\end{align}
\end{mydefn}
\noindent
We refer to $\sups{N}\mathcal{P}_{1+j}(x,y)$ as Genus Two Generalised Weierstrass Functions.
Applying these definitions to \eqref{Zn} we obtain our main theorem:
\begin{theorem}\label{thm:npZhured}
[Quasi-Primary Genus Two Zhu Recursion] The genus two $n$-point  function for a quasi-primary vector $v$ of weight $\wt{[v]}=N$ inserted at $x\in\widehat{\mathcal{S}}_1$ and general vectors $a_1,\ldots,a_L$ and $b_1,\ldots,b_R$  inserted at $x_1,\ldots,x_L\in \widehat{\mathcal{S}}_1$ and $y_1,\ldots,y_R\in \widehat{\mathcal{S}}_2$, respectively, obeys the formal recursive identity
\begin{align}
Z^{(2)}_V (v,x;\boldsymbol{a_l,x_l}|\boldsymbol{b_r,y_r}) 
=&\supsin{N}\mathcal{F}_1(x)\,
O_1(v;\boldsymbol{a_l,x_l}|\boldsymbol{b_r,y_r}) 
\notag
\\
&+ \supsin{N}\mathcal{F}_2(x)\,
O_2(v;\boldsymbol{a_l,x_l}|\boldsymbol{b_r,y_r})
\notag
 \\
&+\supsin{N}\mathcal{F}^\Pi(x)\,
\mathbb{X}_1^\Pi(v;\boldsymbol{a_l,x_l}|\boldsymbol{b_r,y_r}) 
\notag
\\
&+ \sum_{l=1}^L \sum_{j\geq0} \sups{N}\mathcal{P}_{1+j}(x,x_l)
Z^{(2)}_V(\ldots;v[j]a_l,x_l;\ldots) 
\notag
\\
&+ \sum_{r=1}^R \sum_{j\geq0} \sups{N}\mathcal{P}_{1+j}(x,y_r)
Z^{(2)}_V(\ldots;v[j]b_r,y_r;\ldots),
\label{eq:npZhured}
\end{align}
for  $O_a$ of  \eqref{eq:Oadef} and $\mathbb{X}_a^\Pi$ of \eqref{eq:Xadef} and \eqref{nX1GamX22}. There is a similar expression for $v$ inserted on $x\in \widehat{\mathcal{S}}_2$. 
\end{theorem}
\begin{remark}\label{rem:ZhuTh1}
There is a clear analogy between the structure of \eqref{eq:npZhured} and original genus one Zhu recursion \eqref{eq:Zhurec} with elliptic Weierstrass functions replaced by genus two generalised Weierstrass functions.
\end{remark}
\begin{remark}\label{rem:ZhuTh2}
The formal coefficient function 
$\supsin{N}\mathcal{F}(x)$ depends on $N=\wt[v]\ge 1$ and the insertion parameter $x$ but is otherwise universal. 
In particular, these terms determine the $x$ dependence of the genus two 1-point function $Z^{(2)}_V (v,x)$. This is in contrast to genus one 1-point functions which are independent of the torus insertion parameter.
Likewise the formal generalised Weierstrass functions depend on $N$ and insertion points but are otherwise universal.
\end{remark}
\begin{remark}\label{rem:ZhuTh3}
We  show in Section~4 that $\supsin{N}\mathcal{F}_a(x)dx^N$ and $\supsin{N}\mathcal{F}^\Pi(x;m)dx^N$, for $m=1,\ldots ,K-1=2N-3$, provide a basis of holomorphic $N$--differentials in the cases $N=1,2$.  
Since $\supsin{1}\mathcal{F}^\Pi(x)=0$ there are two such terms for $N=1$  whereas for $N\ge 2$ there are $2N-1$ such terms.
This counting agrees with the dimension of the space of genus two holomorphic $N$--differentials following the  Riemann-Roch theorem   \cite{FK}.  
We also show in Sections~4  and  5 that $\sups{N}\mathcal{P}_{1+j}(x,y)$ is holomorphic for $x \neq y$ on the sewing domain in the cases $N=1,2$. The case $N=2$ is particularly significant since this leads to genus two Ward identities with analytic coefficients explored further in Sections~6  and  7.  
\end{remark}
\begin{remark}\label{rem:ZhuTh4}
Following Remark~\ref{rem:M1M2}, there is a corresponding formal genus two Zhu reduction formula to
\eqref{eq:npZhured} for any pair of $V$--modules $M_1,M_2$  involving precisely  the same universal $\supsin{N}\mathcal{F}(x)$ and  $\sups{N}\mathcal{P}_{1+j}(x,\cdot)$  terms.
\end{remark}
\medskip

We may also discuss Zhu reduction for a level $i$   descendant $v^i= \frac{(-1)^i}{i!}L[-1]^iv$ of a quasiprimary vector $v$ inserted at $x\in \widehat{\mathcal{S}}_1$. 
Using translation \eqref{eq:YT} we find
\begin{align}
Z^{(2)}_V (v^i ,x;\boldsymbol{a_l,x_l}|\boldsymbol{b_r,y_r})
=\frac{(-1)^i}{i!}\partial_x^iZ^{(2)}_V (v,x;\boldsymbol{a_l,x_l}|\boldsymbol{b_r,y_r}).
\label{eq:nptqprim}
\end{align} 
The right hand side of \eqref{eq:nptqprim} can be more explicitly expressed in the following way.
Define for $i\ge 0, j>0$ the derivative functions
\begin{align}
\sups{N}\mathcal{P}_{i,1+j}(x,y)&=\frac{(-1)^i}{(i+j)!}\partial_x^i\partial_y^j\left(\sups{N}\mathcal{P}_{1}(x,y)\right)
\notag
\\
&=
\begin{cases}
P_{1+i+j}(x-y)
+ \dfrac{j!(-1)^{1+i+j }}{(i+j)!}\partial_x^i\sups{N}\mathbb{Q}(x)\widetilde{\Lambda}_{\abar} \,\mathbb{P}_{1+j}  (y),\;\mbox{for }x,y\in\widehat{\mathcal{S}}_a,\\
\dfrac{j!(-1)^{N+i+j+1}}{(i+j)!}\partial_x^i\sups{N}\mathbb{Q}(x) \mathbb{P}_{1+j}  (y)
,\;\mbox{for }x\in\widehat{\mathcal{S}}_a,\, y\in\widehat{\mathcal{S}}_{\abar}.
\end{cases}
\label{eq:P2i1j}
\end{align}
Since $Y[v^i,z]   =\frac{ (-1)^i}{i!}\partial_z^iY[v,z]$ it follows that 
\begin{align}
v^i[i+j] = \binom{i+j}{i}v[j],
\label{eq:vvimodes}
\end{align}
for all $j\ge 0$. Hence we find:
\begin{corollary}\label{cor:npointrecursion}
[General Genus Two Zhu Recursion] The genus two $n$-point  function for a level $i\ge 0$ descendant $\frac{(-1)^i}{i!}L[-1]^iv$ of a quasi-primary vector $v$ of weight $\wt{[v]}=N$ inserted at $x\in\widehat{\mathcal{S}}_1$ and general vectors $a_1,\ldots,a_L$ and $b_1,\ldots,b_R$  inserted at $x_1,\ldots,x_L\in \widehat{\mathcal{S}}_1$ and $y_1,\ldots,y_R\in \widehat{\mathcal{S}}_2$, respectively, obeys the recursive identity
\begin{align}
Z^{(2)}_V &\left(\frac{(-1)^i}{i!}L[-1]^iv,x;\boldsymbol{a_l,x_l}|\boldsymbol{b_r,y_r}\right)\notag =
\\
 &\quad\frac{ (-1)^i}{i!}\partial_x^i
\left(\supsin{N}\mathcal{F}_1(x)\right)
O_1(v;\boldsymbol{a_l,x_l}|\boldsymbol{b_r,y_r}) 
\notag
\\
&+\frac{ (-1)^i}{i!}\partial_x^i\left(\supsin{N}\mathcal{F}_2(x)\right)
\,
O_2(v;\boldsymbol{a_l,x_l}|\boldsymbol{b_r,y_r})
\notag
 \\
&+\frac{ (-1)^i}{i!}\partial_x^i\left(\supsin{N}\mathcal{F}^\Pi (x)\right)
\mathbb{X}_1^\Pi(v;\boldsymbol{a_l,x_l}|\boldsymbol{b_r,y_r}) 
\notag
\\
&+ \sum_{l=1}^L \sum_{j\geq0} \ \sups{N}\mathcal{P}_{i,1+j}(x,x_l)
Z^{(2)}_V(\ldots;v[i+j]a_l,x_l;\ldots) 
\notag
\\
&+ \sum_{r=1}^R \sum_{j\geq0} \ \sups{N}\mathcal{P}_{i,1+j}(x,y_r)
Z^{(2)}_V(\ldots;v[i+j]b_r,y_r;\ldots),
\label{eq:npZhuredGen}
\end{align}
for  $O_a$ of  \eqref{eq:Oadef} and $\mathbb{X}_a^\Pi$ of \eqref{eq:Xadef} and \eqref{nX1GamX22}. A similar expression holds for $\frac{(-1)^i}{i!}L[-1]^iv$  inserted at $x\in \widehat{\mathcal{S}}_2$.
\end{corollary}


\section{Holomorphic Weight One Genus Two Zhu Reduction}
\label{sec:Zhuwt1}
\subsection{Identifying the $\supsin{1}\mathcal{F}_a(x)$ and $\sups{1}\mathcal{P}_{1}(x,y)$ coefficients}
In this section we specialize Theorem~\ref{thm:npZhured} to the case where $v$ is a quasi-primary of weight $\wt[v]=N=1$. This implies that  $\Gamma, \Delta, \Pi$ of \eqref{eq:GamDelTh} are given by
\begin{align*}
\Gamma=\Pi=0,\quad \Delta= \mathbbm{1}.
\end{align*}
We can relate the $\supsin{1}\mathcal{F}_a(x)$ coefficients in the  genus two Zhu reduction formula \eqref{eq:npZhured} to the  holomorphic 1-differentials $\nu_a(x)$ and the genus two generalised Weierstrass function $\sups{1}\mathcal{P}_{1}(x,y)$  to the normalised differential of the second kind $\omega(x,y)$ described in Theorem~\ref{Theorem_om_eps} as follows.

Recall  the 1--differentials $a(x)$ of \eqref{eq:adef} and use  \eqref{eq:Sdef} and \eqref{eq:Rdef} to find
\begin{align*}
a(x)=\mathbb{R}(x) S \,dx .
\end{align*}
Hence it follows from
\eqref{eq:Qdef} that
\begin{align*}
\sups{1\,}\mathbb{Q}(x) dx
=
a(x) \left( \mathbbm{1} - A_{\overline{a}}A_a \right)^{-1}S^{-1},
\end{align*}
   for $x\in \widehat{\mathcal{S}}_a$ (recalling that $\Delta= \mathbbm{1}$).
Thus Theorem~\ref{Theorem_om_eps} and \eqref{eq:calFadef} imply
\begin{proposition}
\label{prop:F1nu}
The $N=1$ genus two Zhu reduction coefficient $\supsin{1}\mathcal{F}_a$ for $a=1,2$  is given by
\begin{align}
\nu_a(x)=\supsin{1}\mathcal{F}_a(x)dx,
\label{eq:Fnu}
\end{align}
for normalised holomorphic 1--differentials $\nu_1,\nu_2$ so that
$\supsin{1}\mathcal{F}_a(x)$ is holomorphic for  $x\in \widehat{\mathcal{S}}_b$ and for $(\tau_1,\tau_2,\epsilon)\in 
\mathcal{D}_{\mathrm{sew}}$.
\end{proposition}
In a similar we can we can identify the generalised Weierstrass term $\sups{1}\mathcal{P}_{2}(x,y)$ of \eqref{calP21+j} with $\omega(x,y)$ as expressed in Theorem~\ref{Theorem_om_eps} to find
\begin{align}
\sups{1}\mathcal{P}_{2}(x,y)\, dx dy=\omega(x,y).
\label{eq:P2om}
\end{align}
On integrating, this implies that
\begin{proposition}
\label{prop:P1om}
The $N=1$ genus two generalised Weierstrass function $\sups{1}\mathcal{P}_{1}(x,y)$ is given by the meromorphic  1--differential
\begin{align}
\sups{1}\mathcal{P}_{1}(x,y)\, dx=\int^{y}\omega(x,\cdot),
\label{eq:P1om}
\end{align}
so that
$\sups{1}\mathcal{P}_{1}(x,y)$ is holomorphic for  $x\in \widehat{\mathcal{S}}_a, y\in \widehat{\mathcal{S}}_b$ with $x\neq y$ and for  $(\tau_1,\tau_2,\epsilon)\in 
\mathcal{D}_{\mathrm{sew}}$.
\end{proposition}
\begin{remark}
\label{rem:1Pij}
Since $\sups{1}\mathcal{P}_{i,1+j}(x,y)=\frac{(-1)^i}{(i+j)!}\partial_x^i\partial_y^j\left(\sups{1}\mathcal{P}_{1}(x,y)\right)$ for all $i,j\ge 0$, all $N=1$ genus two generalised Weierstrass functions are holomorphic for all $x\in \widehat{\mathcal{S}}_a, y\in \widehat{\mathcal{S}}_b$ with $x\neq y$ and for all $(\tau_1,\tau_2,\epsilon)\in 
\mathcal{D}_{\mathrm{sew}}$. 
\end{remark}


\subsection{Genus two Heisenberg $n$-point functions}
Consider the rank 1 Heisenberg VOA $M$ generated by $h$ with commutator
\begin{align*}
[h(m),h(n)] = m\delta_{m,-n}.
\end{align*}
The genus two partition function  (found by combinatorial methods) is \cite{MT3}
\begin{align}
Z_{M}^{(2)}(\tau _{1},\tau _{2},\epsilon ) = \frac{1}{\eta(\tau_1)\eta(\tau_2)}\det\left( \mathbbm{1} - A_{1}A_2 \right)^{-1/2},
\label{eq:Z2M}
\end{align}
where $Z_{M}^{(1)}(\tau  ) = \eta(\tau)^{-1}$ for Dedekind eta function 
$\eta(\tau)=q^{1/24}\prod_{n\ge 1}(1-q^n)$.
Hence, by Theorem~\ref{Theorem_A1A2}, $Z_{M}^{(2)}$ is holomorphic on $\mathcal{D}_{\mathrm{sew}}$.
For a pair of irreducible $M$--modules $M_{\lambda_1}=M\otimes e^{\lambda_1}$ and $M_{\lambda_2}=M\otimes e^{\lambda_2}$ one finds  the partition function (cf. Remark~\ref{rem:ZhuTh4}) is given by \cite{MT3}
\begin{align}
Z_{\boldsymbol{\lambda}}^{(2)}(\tau _{1},\tau _{2},\epsilon ) 
:&= \sum_{u\in M} Z_{M_{\lambda_1}}^{(1)}(u;\tau_{1})Z_{M_{\lambda_2}}^{(1)}({\overline u};\tau_{2}) \nonumber \\ 
&= e^{i\pi \boldsymbol{\lambda}\cdot\Omega \cdot\boldsymbol{\lambda}} Z_{M}^{(2)}(\tau _{1},\tau _{2},\epsilon ),
\label{eq:ZMalpha}
\end{align}
where $\boldsymbol{\lambda}\cdot\Omega \cdot\boldsymbol{\lambda}=\sum_{a,b=1}^2\lambda_a\Omega_{ab}\lambda_b$ for the genus two period matrix $\Omega$.
We compute the $1$-point function for the Heisenberg generator $h$ 
\begin{align*}
Z_{\boldsymbol{\lambda}}^{(2)}(h,x;\tau _{1},\tau _{2},\epsilon ) = 
\sum_{u\in M} Z_{M_{\lambda_1}}^{(1)}(Y[h,x]u;\tau_{1})Z_{M_{\lambda_2}}^{(1)}({\overline u};\tau_{2}),
\end{align*}
by means of  Theorem~\ref{thm:npZhured}. We first note that in this case
\begin{align*}
O_1(h;\tau_1,\tau_2,\epsilon) &= \sum_{u\in M} \Tr_{{M_{\lambda_1}}} \left( o(h)o(u) q_1^{L(0) - c/24} \right) Z^{(1)}_{M_{\lambda_2}}(\overline{u};\tau_2) \\ 
&= \lambda_1 \sum_{u\in M} Z_{M_{\lambda_1}}^{(1)}(u;\tau_{1})Z_{M_{\lambda_2}}^{(1)}({\overline u};\tau_{2}) \\
&= \lambda_1 Z_{\boldsymbol{\lambda}}^{(2)}(\tau _{1},\tau _{2},\epsilon ),
\end{align*}
and similarly $O_2(h;\tau_1,\tau_2,\epsilon) = \lambda_2 Z_{\boldsymbol{\lambda}}^{(2)}(\tau _{1},\tau _{2},\epsilon )$. 
Defining 
\begin{align}
\nu_{\boldsymbol{\lambda}}(x)=\lambda_1 \nu_1(x) + \lambda_2 \nu_2(x),
\label{eq:nualpha}
\end{align}
and applying  \eqref{eq:Fnu} to Theorem~\ref{thm:npZhured}  we obtain
\begin{proposition}\label{HeisenbergM1point}
The $1$-point correlation function for the Heisenberg generator $h$ for a pair of irreducible Heisenberg modules $M_{\lambda_1},M_{\lambda_2}$ is given by
\begin{align*}
Z_{\boldsymbol{\lambda}}^{(2)}(h,x;\tau _{1},\tau _{2},\epsilon ) dx 
=\nu_{\boldsymbol{\lambda}}(x)Z_{\boldsymbol{\lambda}}^{(2)}(\tau _{1},\tau _{2},\epsilon ).
\end{align*}
\end{proposition}
\noindent
This  agrees  with Theorem~12 of \cite{MT3} obtained by a combinatorial method.

\medskip

We next consider, for two irreducible modules $M_{\lambda_1},M_{\lambda_2}$, the $n$-point function for $L+1$ Heisenberg vectors $h$  inserted at $x,x_1,\ldots,x_{L}\in \widehat{\mathcal{S}}_1$ and $R$  
Heisenberg vectors $h$  inserted at $y_1,\ldots,x_{R}\in \widehat{\mathcal{S}}_2$
\begin{align*}
Z_{\boldsymbol{\lambda}}^{(2)} (h,x;\boldsymbol{h,x_l}|\boldsymbol{h,y_r})
= \sum_{u\in M} Z^{(1)}_{M_{\lambda_1}}(h,x;\boldsymbol{Y[h,x_l]}u;\tau_1)Z^{(1)}_{M_{\lambda_2}}(\boldsymbol{Y[h,y_r]}\overline{u};\tau_2).
\end{align*}
Much as before we find that
\begin{align*}
O_a(h;\boldsymbol{h,x_l}|\boldsymbol{h,y_r}) 
= \lambda_a Z_{\boldsymbol{\lambda}}^{(2)} (\boldsymbol{h,x_l}|\boldsymbol{h,y_r}),\quad a=1,2,
\end{align*}
for Heisenberg $(L+R)$--point function $Z_{\boldsymbol{\lambda}}^{(2)} (\boldsymbol{h,x_l}|\boldsymbol{h,y_r})$.
Furthermore, since $h[j]h = \delta_{j1} \mathbf{1}$ and using \eqref{eq:P2om} 
then Theorem~\ref{thm:npZhured} implies:
\begin{proposition}\label{HeisenbergMnpoint}
The genus two $n$-point function for $L+R+1$ Heisenberg vectors $h$ inserted at $x,x_1,\ldots,x_{n-1}\in \widehat{\mathcal{S}}_1$ and at $y_1,\ldots,x_{R}\in \widehat{\mathcal{S}}_2$ and for irreducible Heisenberg modules $M_{\lambda_1},M_{\lambda_2}$ is given by
\begin{align*}
&Z_{\boldsymbol{\lambda}}^{(2)} (h,x;\boldsymbol{h,x_l}|\boldsymbol{h,y_r})
dx \prod_{k=1}^{L}dx_k  \prod_{s=1}^{R}dy_s \\ 
&=\nu_{\boldsymbol{\lambda}}(x)
Z_{\boldsymbol{\lambda}}^{(2)} (\boldsymbol{h,x_l}|\boldsymbol{h,y_r})
\prod_{k=1}^{L}dx_k  \prod_{s=1}^{R}dy_s \\
&\; + \sum_{l=1}^{L} \omega(x,x_l) 
Z_{\boldsymbol{\lambda}}^{(2)} (h,x_1;\ldots;\widehat{h,x_l};\ldots;h,x_L|\boldsymbol{h,y_r}) 
\prod_{k=1,k\neq l }^{L}dx_k  \prod_{s=1}^{R}dy_s\\
&\; + \sum_{r=1}^{R} \omega(x,y_r) 
Z_{\boldsymbol{\lambda}}^{(2)} (\boldsymbol{h,x_l}|h,y_1;\ldots;\widehat{h,y_r};\ldots;h,y_R) 
\prod_{k=1 }^{L}dx_k  \prod_{s=1,s\neq r}^{R}dy_s,
\end{align*}
where $\widehat{h,x}$ etc. denotes omission of the given term.
\end{proposition}
\noindent
This  agrees with Theorem 13 of \cite{MT3} proved by a combinatorial method. 
In fact, all $n$--point functions for the Heisenberg VOA are generated by such $n$--point functions for the Heisenberg vector $h$ \cite{MT3}.


\section{Weight Two Genus Two Zhu Reduction}
\label{sec:Zhuwt2}

In this section we specialise Theorem~\ref{thm:npZhured} to the case where $v$ is quasi-primary of weight $\wt[v]=N=2$. We demonstrate the holomorphy of the $\supsin{2}\mathcal{F}(x)$ terms which appear. This allows us to express genus two Ward identities and genus two Virasoro $n$-point correlation functions in terms of a covariant derivative with respect to the parameters $(\tau_1,\tau_2,\epsilon)$ and the generalised Weierstrass functions. 
The general expressions derived here are examined further in Sect.~\ref{sec:GDE}, where analysis of the resulting partial differential equations demonstrates the holomorphy of all coefficient terms appearing in $N=2$ genus two Zhu reduction.     

\subsection{Holomorphy of $\supsin{2}\mathcal{F}(x)$ terms}

We now specialize Theorem~\ref{thm:npZhured} to the case where $v$ is quasi-primary of weight $\wt[v]=N=2$ so that the infinite matrices of $\Gamma, \Delta, \Pi$ of \eqref{eq:GamDelTh} are 
\begin{align*}
\Gamma=\Pi=
\begin{bmatrix}
1 & 0 & 0&\ldots\\
0 &0 & 0 &\ldots\\
0 &0 & 0 &\ldots\\
\vdots & \vdots & \vdots& \ddots
\end{bmatrix}
,\quad \Delta= \begin{bmatrix}
0 & 0 & 0 & 0 &\ldots\\
0 & 0 & 0 & 0 &\ldots\\
1 & 0 & 0 & 0 &\ldots\\
0 & 1 & 0 & 0  &\ldots\\
0 & 0 & 1 & 0  &\ldots\\
\vdots & \vdots & \vdots & \vdots &\ddots
\end{bmatrix}.
\end{align*}
We now relate the formal $\supsin{2}\mathcal{F}_a(x)$ and $\supsin{2}\mathcal{F}^\Pi(x)$ coefficients appearing in the $N=2$ genus two Zhu reduction formula \eqref{eq:npZhured}
to the three dimensional  space of holomorphic 2--differentials. 
The geometric meaning of the genus two Weierstrass function $\sups{2\,}\mathcal{P}_{1}(x,y)$ will be described later on  in Proposition~\ref{prop:P21}.

Let $\{\Phi_r(x)\}$, for $r=1,2,3$, denote the formal 2--differentials
\begin{align}
\Phi_1(x)=\supsin{2}\mathcal{F}_1(x)dx^2,\quad 
\Phi_2(x)=\supsin{2}\mathcal{F}_2(x)dx^2,\quad 
\Phi_3(x)=\epsilon^{-1/2}\cdot\supsin{2}\mathcal{F}^\Pi(x;1)dx^2.
\label{eq:Phibasis}
\end{align}
\begin{theorem}
\label{thm:Psi2forms}
$\{\Phi_r(x)\}$ is a  basis of holomorphic 2--differentials  for $x\in \widehat{\mathcal{S}}_a$ and for $(\tau_1,\tau_2,\epsilon)\in \mathcal{D}_{\mathrm{sew}}$ with normalization
\begin{align}
\frac{1}{\tpi}\oint_{\alpha^i} \Phi_r(z) (dz)^{-1} &=\delta_{ri},
\quad \frac{1}{\tpi}\oint_{\mathcal{C}_a}z_a\, \Phi_r(z_a) (dz_a)^{-1} =\delta_{r3},
\label{eq:PsiNorm}
\end{align}
with $i=1,2$ 
where $\alpha^i$ is the standard genus two homology cycle and $\mathcal{C}_a $
is an anti-clockwise contour surrounding the excised disc centred at $z_a=0$ on  $\widehat{\mathcal{S}}_a$. 
\end{theorem}
\begin{remark}
\label{rem:normPhi}
The normalization conditions \eqref{eq:PsiNorm} are all coordinate dependent unlike the analogous  condition \eqref{omeganorm} for holomorphic 1--differentials.
\end{remark}

\noindent\textbf{{Proof.}}
Let $\Psi(x)=(\Psi_r(x))$, for $ r=1,2,3$, denote the column vector with components given by the 3 independent genus two holomorphic 2--differentials
\begin{align}
\Psi_1(x)=\nu_1(x)^2,\quad 
\Psi_2(x)=\nu_2(x)^2,\quad 
\Psi_3(x)=\nu_1(x)\nu_2(x),
\label{eq:Psibasis}
\end{align}
for normalised 1--differentials $\nu_i(x)$. Let $\Xi=(\Xi_{rs})$ denote the  2--differential period matrix over the cycles $\alpha^i$ and the contour $\mathcal{C}_a$ defined by
\begin{align}
\Xi_{ri}&=
\frac{1}{\tpi}\oint_{\alpha^i} \Psi_r(z) (dz)^{-1}\quad \mbox{for } i=1,2,
\notag
\\
\Xi_{r3}&=
\frac{1}{\tpi}\oint_{\mathcal{C}_a} z_a\, \Psi_r(z_a) (dz_a)^{-1}.
\label{eq:Xi}
\end{align}
Clearly, $\Xi_{rs}$ is  holomorphic in $(\tau_1,\tau_2,\epsilon)\in \mathcal{D}_{\mathrm{sew}}$. 
With $z_2=\epsilon/z_1$ we note that  
\begin{align*}
\frac{1}{\tpi}\oint_{\mathcal{C}_1} z_1\, \Psi_r(z_1) (dz_1)^{-1}
=&
\frac{1}{\tpi}\left(-\oint_{\mathcal{C}_2} \right)
 \left(\frac{\epsilon}{z_2}\right) \Psi_r(z_2)
 \left(-\frac{\epsilon}{z_2^2}dz_2\right)^{-1}\\
=&
\frac{1}{\tpi}\oint_{\mathcal{C}_2}
z_2\, \Psi_r(z_2)
 \left(dz_2\right)^{-1},
\end{align*}
so that $\Xi_{r3}$ is well--defined.
We now show  that $\Xi$ is invertible on $\mathcal{D}_{\mathrm{sew}}$ and that
\begin{align}
\Phi(x)=\Xi^{-1}\Psi(x),
\label{eq:PhiPsi}
\end{align}
where $\Phi(x)=(\Phi_r(x))$. \eqref{eq:PhiPsi}  implies that $\{\Phi_r(x)\}$ is a indeed basis of holomorphic  2--differentials with  normalization \eqref{eq:PsiNorm}.

In order to prove \eqref{eq:PhiPsi}, we compute in two separate ways the genus two 1--point function for a particular weight $N=2$ primary vector in the rank 2 Heisenberg VOA $M^2=M\otimes M$. 
One computation is manifestly holomorphic and expressed in terms of the  2--differential basis $\{\Psi_r(x)\}$ whereas the other  is in terms of the formal  2--differentials $\{\Phi_r(x)\}$.
$M^2$ is generated by $h^1=h\otimes \vac$ and $h^2=\vac \otimes h$. 
Let $M_{\lambda_1,\mu_1}$ and $M_{\lambda_2,\mu_2}$ (where $M_{\lambda,\mu}=M_{\lambda}\otimes M_{\mu}$) be a pair of irreducible modules for this VOA.  
Using the shorthand notation 
\begin{align*}
Z^{(2)}_{\boldsymbol{\lambda},\boldsymbol{\mu}} (\ldots) = Z^{(2)}_{M_{\lambda_1,\mu_1}M_{\lambda_2,\mu_2 }} (\ldots) , 
\end{align*}
we then find that Proposition~\ref{HeisenbergMnpoint} implies
\begin{align*}
Z^{(2)}_{\boldsymbol{\lambda},\boldsymbol{\mu}} (h^1,x;h^2,y)\,dxdy=
\nu_{\boldsymbol{\lambda}}(x)\nu_{\boldsymbol{\mu}}(y)\, Z^{(2)}_{\boldsymbol{\lambda},\boldsymbol{\mu}},
\end{align*}
which is holomorphic  for $x\in \widehat{\mathcal{S}}_a$, $y\in \widehat{\mathcal{S}}_b$ and 
$(\tau_1,\tau_2,\epsilon)\in \mathcal{D}_{\mathrm{sew}}$.
Taking $x=y$ gives the genus two 1--point function for 
$v=h \otimes h$, a weight 2 primary vector with vertex operator $Y(v,x)=Y(h,x)\otimes Y(h,x)$ where 
\begin{align}
Z^{(2)}_{\boldsymbol{\lambda},\boldsymbol{\mu}} (v,x)dx^2=&
\nu_{\boldsymbol{\lambda}}(x)\nu_{\boldsymbol{\mu}}(x)Z^{(2)}_{\boldsymbol{\lambda},\boldsymbol{\mu}}.
\label{eq:Znunu}
\end{align}
Alternatively, Theorem~\ref{thm:npZhured} for $N=2$ and \eqref{eq:Psibasis} implies that
for $x\in \widehat{\mathcal{S}}_1$
\begin{align}
Z^{(2)}_{\boldsymbol{\lambda},\boldsymbol{\mu}} (v,x) dx^2
= \,\Phi_1(x)O_1(v) 
+ \Phi_2(x) O_2(v)
+\Phi_3(x)\cdot \epsilon^{1/2}\mathbb{X}_1^\Pi(v;1).
\label{eq:Zv2}
\end{align}
From \eqref{eq:Oadef} and \eqref{eq:Xadef} we have
\begin{align*}
O_1(v) =&\sum_{u\in M^2} \Tr_{M_{\lambda_1,\mu_1}} \left( o(v)o(u)q_1^{L(0) - c/24} \right) 
Z^{(1)}_{M_{\lambda_2,\mu_2}}(\overline{u},\tau_2),
\\
O_2(v) =&\sum_{u\in M^2}  
Z^{(1)}_{M_{\lambda_2,\mu_2}}(u,\tau_1)
\Tr_{M_{\lambda_1,\mu_1}} \left( o(v)o(\overline{u})q_2^{L(0) - c/24} \right),
\\
\epsilon^{1/2}\mathbb{X}_1^\Pi(v;1)=&\sum_{u\in M^2} Z_{M_{\lambda_1,\mu_1}} 
(v[1]u,\tau_1)
Z^{(1)}_{M_{\lambda_2,\mu_2}}(\overline{u},\tau_2).
\end{align*}
But these terms may be found from the holomorphic expression \eqref{eq:Znunu}. In particular, using $o(v)=v(1)$ (since $v$ is a primary vector of weight 2) we find
\begin{align*}
O_a(v) =&
 \frac{1}{\tpi}\int_{0}^{\tpi}Z^{(2)}_{\boldsymbol{\lambda},\boldsymbol{\mu}} (v,x) dx
\\
=&
Z^{(2)}_{\boldsymbol{\lambda},\boldsymbol{\mu}}\cdot \frac{1}{\tpi}\oint_{\alpha^a}\nu_{\boldsymbol{\lambda}}(x)\nu_{\boldsymbol{\mu}}(x)(dx)^{-1},
\end{align*} 
and 
\begin{align*}
\epsilon^{1/2}\mathbb{X}_1^\Pi(v;1)=&
 \frac{1}{\tpi}\oint_{\mathcal{C}_1}x\,Z^{(2)}_{\boldsymbol{\lambda},\boldsymbol{\mu}} (v,x) dx
\\
=&
Z^{(2)}_{\boldsymbol{\lambda},\boldsymbol{\mu}}\cdot
 \frac{1}{\tpi}
\oint_{\mathcal{C}_1}\nu_{\boldsymbol{\lambda}}(x_1)\nu_{\boldsymbol{\mu}}(x_1)(dx_1)^{-1}.
\end{align*}
Thus comparing the respective $\lambda_1\mu_1, \lambda_2\mu_2$ and $\lambda_1\mu_2+\lambda_2\mu_1$ terms in \eqref{eq:Znunu} and  \eqref{eq:Zv2} we find the holomorphic  2--differentials \eqref{eq:Psibasis} are given by
\begin{align}
\Psi(x)=\Xi\Phi(x),
\label{eq:PhiPsi2}
\end{align}
where $\Xi$ is the holomorphic  2--differential period matrix  \eqref{eq:Xi}. 

We lastly show that $\Xi$ is invertible on $\mathcal{D}_{\mathrm{sew}}$. Suppose that $\Xi$ is singular for some $(\tau_1,\tau_2,\epsilon)\in \mathcal{D}_{\mathrm{sew}}$. 
Then there must exist a row vector $\kappa\neq 0$ for which $\kappa\Xi=0$. Hence \eqref{eq:PhiPsi2} would imply that $\kappa\Psi(x)=0$ which contradicts the fact that $\{\Psi_r(x)\}$ is an independent basis of holomorphic  2--differentials.   Therefore \eqref{eq:PhiPsi} is true and the theorem follows. \qed

\subsection{Genus two Virasoro $1$-point functions}
We next consider applications of Theorems~\ref{thm:npZhured} and \ref{thm:Psi2forms} to the important case where $v=\omt$, the square bracket VOA Virasoro vector of weight $N=\wt[v]=2$.
We first consider the Virasoro vector $1$-point function inserted at 
 $x\in \widehat{\mathcal{S}}_1$  for which  Theorem~\ref{thm:npZhured} for $N=2$ implies (much as in \eqref{eq:Zv2})
\begin{align}
Z^{(2)}_{V} (\omt,x) dx^2
= \,\Phi_1(x)O_1(\omt) 
+ \Phi_2(x) O_2(\omt)
+\Phi_3(x)\cdot \epsilon^{1/2}\mathbb{X}_1^\Pi(\omt;1),
\label{eq:Zom}
\end{align} 
for the holomorphic  2--differentials $\Phi_r(x)$ of \eqref{eq:Phibasis}. Here
\begin{align*} O_1(\omt) 
&=\sum_{u\in V} \Tr_V \left( o(\omt)o(u)q_1^{L(0)-c/24}\right) Z^{(1)}_V(\overline{u};\tau_2) \\
&= \sum_{u\in V}  q_1\partial_{q_1}Z^{(1)}_V(u;\tau_1) Z^{(1)}_V(\overline{u};\tau_2) \\
&= q_1\partial_{q_1} Z^{(2)}_V (\tau_1,\tau_2,\epsilon),
\end{align*}
and similarly  $O_2(\omt) = q_2\partial_{q_2} Z^{(2)}_V (\tau_1,\tau_2,\epsilon)$. 
Furthermore, since $\omt[1]=L[0]$
\begin{align*}
\epsilon^{1/2} \mathbb{X}_1^\Pi(\omt;1) 
&= \sum_{u\in V} Z^{(1)}_V(L[0]u;\tau_1) Z^{(1)}_V(\overline{u};\tau_2).
\end{align*}
Choose an $L[0]$--homogeneous basis $\{ u \}$ giving
\begin{align*}
\epsilon^{1/2}\mathbb{X}_1^\Pi(1) 
&=  \sum_{n\geq0}\sum_{u\in V_{[n]}} nZ^{(1)}_V(u;\tau_1) Z^{(1)}_V(\overline{u};\tau_2).
\end{align*}
$Z^{(1)}_V(u;\tau_1) Z^{(1)}_V(\overline{u};\tau_2)$ is proportional to $\epsilon^n$ from  \eqref{eq:LiZA} for $u\in V_{[n]}$  and  so
\begin{align*}
\epsilon^{1/2}\mathbb{X}_1^\Pi(1) 
&= \sum_{n\geq0}\sum_{u\in V_{[n]}}\epsilon\partial_\epsilon \left(Z^{(1)}_V(u;\tau_1) Z^{(1)}_V(\overline{u};\tau_2) \right)= \epsilon\partial_\epsilon Z^{(2)}_V(\tau_1,\tau_2,\epsilon).
\end{align*}
We therefore define the differential operator
\begin{align}
\boldsymbol{D}_{x} 
= \supsin{2}\mathcal{F}_1(x)\,q_1\partial_{q_1}
 + \supsin{2}\mathcal{F}_2(x) \, q_2\partial_{q_2} 
+   \supsin{2}\mathcal{F}^\Pi(x;1)\,\epsilon^{1/2} \partial_\epsilon,
\label{eq:Dx}
\end{align}
where the subscript  indicates the dependence on $x\in \widehat{\mathcal{S}}_a$.  
Equivalently, 
\begin{align} 
dx^2\boldsymbol{D}_{x}
= \Phi_1(x)\,q_1\partial_{q_1}
 + \Phi_2(x) \, q_2\partial_{q_2} 
+   \Phi_3(x) \,\epsilon \partial_\epsilon,
\label{eq:Dxdx}
\end{align}
for holomorphic  2--differentials $\{\Phi_r(x)\}$ with normalization \eqref{eq:PsiNorm}.
Thus we have
\begin{proposition}\label{Virasoro1pt}
The genus two Virasoro $1$-point correlation function for a VOA $V$ is given by
\begin{align}
Z^{(2)}_V(\omt,x;\tau_1,\tau_2,\epsilon) = \boldsymbol{D}_x Z^{(2)}_V(\tau_1,\tau_2,\epsilon).
\label{eq:Zom2}
\end{align}
\end{proposition}
\noindent
Proposition~\ref{Virasoro1pt} is analogous to the genus one Virasoro 1--point function  \eqref{eq:delqZ}.
Note that $\boldsymbol{D}_x$ acts on differentiable functions on $\mathcal{D}^\epsilon$. 
We will directly relate $\boldsymbol{D}_x$ to the differential operator $\nabla_x$ of \eqref{eq:nabla}  in the next section. 

%


\subsection{Genus two Ward identities}
Consider the $n$-point function for $\omt$ inserted at $x$ and Virasoro primary vectors $a_1,\ldots,a_L$ and $b_1,\ldots,b_R$ (of respective $L[0]$ weight $ \wt[a_1],\ldots , \wt[b_R]$)  inserted at $x_1,\ldots,x_L\in \widehat{\mathcal{S}}_1$  and $y_1,\ldots,y_R\in \widehat{\mathcal{S}}_2$  respectively.  Since $L[j-1]=\omt[j]$ and $a_l,b_r$  are primary vectors, we find Theorem~\ref{thm:npZhured} implies 
\begin{align}
Z^{(2)}_V (\omt,x;\boldsymbol{a_l,x_l}|\boldsymbol{b_r,y_r}) 
&= \sups{2}\mathcal{F}_1(x)\,
O_1(\omt;\boldsymbol{a_l,x_l}|\boldsymbol{b_r,y_r}) 
\notag
\\
&+ \sups{2}\mathcal{F}_2(x)\,
O_2(\omt;\boldsymbol{a_l,x_l}|\boldsymbol{b_r,y_r})
\notag
 \\
&+\sups{2}\mathcal{F}^\Pi(x)\,
\mathbb{X}_1^\Pi(\omt;\boldsymbol{a_l,x_l}|\boldsymbol{b_r,y_r}) 
\notag
\\
&+ \sum_{l=1}^L \sups{2}\mathcal{P}_{1}(x,x_l)
Z^{(2)}_V(\ldots;L[-1]a_l,x_l;\ldots) 
\notag
\\&+ \sum_{l=1}^L \sups{2}\mathcal{P}_{2}(x,x_l)
Z^{(2)}_V(\ldots;L[0]a_l,x_l;\ldots) 
\notag
\\
&+ \sum_{r=1}^R  \sups{2}\mathcal{P}_{1}(x,y_r)
Z^{(2)}_V(\ldots;L[-1]b_r,y_r;\ldots)
\notag
\\
&+ \sum_{r=1}^R  \sups{2}\mathcal{P}_{2}(x,y_r)
Z^{(2)}_V(\ldots;L[0]b_r,y_r;\ldots).
\label{eq:WardZhu}
\end{align}
Much as for the Virasoro 1-point function we find that
\begin{align*}
O_a(\omt;\boldsymbol{a_l,x_l}|\boldsymbol{b_r,y_r}) 
=&q_a\partial_{q_a} Z^{(2)}_V (\boldsymbol{a_l,x_l}|\boldsymbol{b_r,y_r}),\\
\epsilon^{1/2}\mathbb{X}_1^\Pi(\omt;1;\boldsymbol{a_l,x_l}|\boldsymbol{b_r,y_r})
= &\epsilon \partial_\epsilon Z^{(2)}_V(\boldsymbol{a_l,x_l}|\boldsymbol{b_r,y_r}),
\end{align*}
and altogether we obtain
\begin{proposition}\label{prop:WardId}
The $n$-point function obeys the genus two Ward identity
\begin{align}
&Z^{(2)}_V (\omt,x;\boldsymbol{a_l,x_l}|\boldsymbol{b_r,y_r};\tau_1,\tau_2,\epsilon)
\notag
 \\
&= \Bigg( \boldsymbol{D}_x 
+ \sum_{l=1}^L \Big( \sups{2}\mathcal{P}_{1}(x,x_l) \partial_{x_l} 
+ \wt[a_l]\cdot\sups{2}\mathcal{P}_{2}(x,x_l)\Big) 
\notag
\\
&\quad \quad
+ \sum_{r=1}^R \Big( \sups{2}\mathcal{P}_{1}(x,y_r) \partial_{y_r} 
+ \wt[b_r]\cdot \sups{2}\mathcal{P}_{2}(x,y_r)\Big)\Bigg)
Z^{(2)}_V(\boldsymbol{a_l,x_l}|\boldsymbol{b_r,y_r}),
\label{eq:WardId}
\end{align}
where $a_1,\ldots,a_L,b_1,\ldots,b_R\in V$ are primary vectors  of respective $L[0]$ weight $ \wt[a_1],\ldots , \wt[b_R]$ and 
with $\boldsymbol{D}_x$ of \eqref{eq:Dx}.
\end{proposition}
\noindent 
\eqref{eq:WardId}  is analogous to the genus one Ward Identity \eqref{eq:Z1Ward}.
\begin{remark}
\label{rem:coeffcgt} 
In Theorem~\ref{th:Zhured2} below we prove the convergence of all the coefficients appearing in \eqref{eq:WardId} for all $x\neq x_l,y_r$ and for all $(\tau_1,\tau_2,\epsilon)\in \mathcal{D}_{\mathrm{sew}}$.
\end{remark}

\medskip 
We lastly consider the Virasoro $n$-point correlation function for $\omt$
inserted at $x$, $x_1,\ldots,x_L \in \widehat{\mathcal{S}}_1$
 and $y_1,\ldots,y_R\in \widehat{\mathcal{S}}_2$  respectively.  We find, much as in 
Proposition~\ref{prop:WardId} and using $L[2]\omt=\frac{c}{2}\vac$, that Theorem~\ref{thm:npZhured} implies
\begin{proposition}\label{prop:Vnpt}
The genus two Virasoro $n$-point correlation  function is
\begin{align}
&Z^{(2)}_V (\omt,x;\boldsymbol{\omt,x_l}|\boldsymbol{\omt,y_r};\tau_1,\tau_2,\epsilon)
\notag
 \\
=& \Bigg( \boldsymbol{D}_x 
+ \sum_{l=1}^L \left( \sups{2}\mathcal{P}_{1}(x,x_l) \partial_{x_l} + 2\cdot\sups{2}\mathcal{P}_{2}(x,x_l) \right) 
\notag
 \\
& \qquad + \sum_{r=1}^R \left( \sups{2}\mathcal{P}_{1}(x,y_r) \partial_{y_r} + 2\cdot\sups{2}\mathcal{P}_{2}(x,y_r) \right) \Bigg) Z^{(2)}_V (\boldsymbol{\omt,x_l}|\boldsymbol{\omt,y_r}) 
\notag
\\
&+ \frac{c}{2} \sum_{l=1}^L \sups{2}\mathcal{P}_{4}(x,x_l) Z^{(2)}_V\left(\omt,x_1;\ldots;\widehat{\omt,x_l};\ldots;\omt,x_L|\boldsymbol{\omt,y_r}\right) 
\notag
\\
&+ \frac{c}{2} \sum_{r=1}^R \sups{2}\mathcal{P}_{4}(x,y_r) Z^{(2)}_V\left(\boldsymbol{\omt,x_l}|\omt,y_1;\ldots;\widehat{\omt,y_r};\ldots;\omt,y_R\right),
\label{eq:Vnpt}
\end{align}
where $\widehat{\omt,x_l}$ and $\widehat{\omt,y_r}$ denotes omission of the given term.
\end{proposition}

\noindent 
Remark~\ref{rem:coeffcgt} again applies  concerning the convergence of the coefficients in \eqref{eq:Vnpt}.
Also note that \eqref{eq:Vnpt}  is analogous to the genus one Ward Identity \eqref{eq:Z1Ward2}.


\section{Analytic Genus Two Differential Equations}
\label{sec:GDE}
In this section, we discuss the geometric significance of the differential operator $\mathbf{D}_x$  \eqref{eq:Dx} which was employed in the  previous section. 
We derive a number of differential equations involving $\mathbf{D}_x$ arising from the Heisenberg VOA.  
One important consequence is a proof that the holomorphic map $F^\epsilon$ from the sewing domain $
\mathcal{D}_{\mathrm{sew}}$ to the Siegel upper half plane $\HH_{2}$ provided by $\Omega (\tau _{1},\tau_{2},\epsilon )  $ is injective  but not surjective and we show that the differential operators $\mathbf{D}_x$ and $\nabla_x$ of \eqref{eq:nabla} are equivalent on the sewing domain. 
Further, the holomorphy of all coefficient terms appearing in the genus two Ward identities and Virasoro $n$-point functions derived in Sect.~\ref{sec:Zhuwt2} is demonstrated. As a direct consequence, the genus two differential equations arising from Virasoro singular vectors have holomorphic coefficients.

As in previous sections, we suppress the dependence on $\tau_1,\tau_2, \epsilon$ where there is no ambiguity.

\subsection{The  injectivity and non-surjectivity of $F^\epsilon$}
Consider the Heisenberg VOA $M$ generated by $h$. Proposition~\ref{HeisenbergMnpoint} implies
\begin{align*}
Z^{(2)}_M(h,x;h,y)\,dxdy=\omega(x,y)Z^{(2)}_M.
\end{align*}
Since $\omt=\frac{1}{2}h[-1]^2\vac$   the Virasoro $1$-point function can be obtained in the limit
\begin{align*}
dx^2 Z^{(2)}_M(\widetilde{\omega},x)
&= \lim_{x\rightarrow y} \frac{1}{2} \left( Z^{(2)}_M(h,x;h,y) - \frac{1}{(x-y)^2} Z^{(2)}_M  \right)dx dy = \frac{1}{12} s(x) Z^{(2)}_M ,
\end{align*}
where $s(x)$ is the projective connection \eqref{eq:projcon}. Comparing with \eqref{eq:Zom2} we find
\begin{proposition}\label{prop:HeisDE}
The genus two partition function for the rank $1$ Heisenberg  VOA satisfies the differential equation
\begin{align}
dx^2 \, \mathbf{D}_x Z^{(2)}_M (\tau_1,\tau_2,\epsilon) 
= \frac{1}{12} s(x) Z^{(2)}_M (\tau_1,\tau_2,\epsilon).
\label{eq:HeisDE}
\end{align}
\end{proposition}
\noindent This  is analogous to the genus one  result $q\partial_q Z^{(1)}_M (\tau) = \frac{1}{2}E_2(\tau) Z^{(1)}_M (\tau)$ for $Z^{(1)}_M (\tau)=1/\eta(q)$ with the projective connection playing the role of $E_2(\tau)$.

\medskip
Similarly, consider the genus two Virasoro $1$--point function for a pair of  Heisenberg modules $M_{\lambda_1},M_{\lambda_2}$ with 
$Z^{(2)}_{\boldsymbol{\lambda}} =e^{i\pi \boldsymbol{\lambda}\cdot\Omega \cdot\boldsymbol{\lambda}} Z^{(2)}_M$. 
In this case, Proposition~\ref{HeisenbergMnpoint} implies
\begin{align*}
Z^{(2)}_{\boldsymbol{\lambda}}(h,x;h,y)\,dxdy=
\left(\nu_{\boldsymbol{\lambda}}(x)\nu_{\boldsymbol{\lambda}}(y)+\omega(x,y)\right)Z^{(2)}_{\boldsymbol{\lambda}}.
\end{align*}
Taking the $x\rightarrow y$ limit  we  find
\begin{align}
 Z^{(2)}_{\boldsymbol{\lambda}} (\widetilde{\omega},x) dx^2
= &\left( \frac{1}{2} \nu_{\boldsymbol{\lambda}}(x)^2 +\frac{1}{12} s(x) 
 \right)
Z^{(2)}_{\boldsymbol{\lambda}}.
\label{eq:Zlamom}
\end{align}
But using   \eqref{eq:ZMalpha} and   \eqref{eq:Zom2} we also  have
\begin{align*}
 Z^{(2)}_{\boldsymbol{\lambda}} (\widetilde{\omega},x) dx^2
&= dx^2 \, \mathbf{D}_x \left( e^{i\pi \boldsymbol{\lambda}\cdot\Omega \cdot\boldsymbol{\lambda}} Z^{(2)}_M  \right) \\
&= \left(  dx^2 \, \mathbf{D}_x \left( i\pi\boldsymbol{\lambda}\cdot\Omega\cdot\boldsymbol{\lambda} \right) 
+\frac{1}{12} s(x) \right) Z^{(2)}_{\boldsymbol{\lambda}}.
\end{align*}
Comparing these expressions we obtain:
\begin{proposition}\label{prop:DxOmega}
For $i,j=1,2$
\begin{align}
2\pi i \, dx^2 \, \mathbf{D}_x \Omega_{ij}  = \nu_i(x)\nu_j(x) .
\label{eq:DxOmega}
\end{align}
\end{proposition}
It is convenient to define $\Omega_{r}$ for  $r=1,2,3$, and $\tau_3$ by
\begin{align*}
\Omega_1=\Omega_{11},\quad \Omega_2=\Omega_{22},\quad \Omega_3=\Omega_{12},\quad \epsilon=e^{\tpi\tau_3}.
\end{align*}
Recall the differential operator  $\nabla_x$ of \eqref{eq:nabla} which we can rewrite as
\begin{align*}
\nabla_x =\frac{1}{2\pi i} \sum_{r=1}^3 \Psi_r(x) \frac{\partial}{\partial \Omega_{r}},
\end{align*}
using the holomorphic 2--differential basis $\{\Psi_r(x)\}$ of \eqref{eq:Psibasis}.  
From \eqref{eq:Dxdx} we may rewrite \eqref{eq:DxOmega} in terms of the bases $\{\Psi_r(x)\}$ and $\{\Phi_r(x)\}$  of \eqref{eq:Phibasis} to find
\begin{align*}
\Psi(x)= \frac{\partial( \Omega_{1},\Omega_{2},\Omega_{3})}{\partial(\tau_1,\tau_2,\tau_3)}\Phi(x),
\end{align*}
for Jacobian matrix $\frac{\partial( \Omega_{1},\Omega_{2},\Omega_{3})}{\partial(\tau_1,\tau_2,\tau_3)}$ and column vectors $\Psi(x)=(\Psi_r(x))$ and $\Phi(x)=(\Phi_r(x))$. 
Referring to  \eqref{eq:PhiPsi2}, this implies
\begin{align*}
\frac{\partial( \Omega_{1},\Omega_{2},\Omega_{3})}{\partial(\tau_1,\tau_2,\tau_3)}=\Xi,
\end{align*} 
where $\Xi$ is  the holomorphic  2--differential period matrix defined in \eqref{eq:Xi}.  
But by Theorem~\ref{thm:Psi2forms},  $\Xi$ is invertible on the sewing domain and therefore, by the inverse function theorem, the map $(\tau _{1},\tau _{2},\epsilon ) \mapsto \Omega (\tau _{1},\tau_2,\epsilon)$ is one to one. Hence, altogether we have the following  result:
\begin{theorem}\label{thm:Dxnabla1}
The holomorphic map 
\begin{align*}
F^{\Omega }:\mathcal{D}_{\mathrm{sew}} &\rightarrow \HH_{2}, 
\\
(\tau _{1},\tau _{2},\epsilon ) &\mapsto \Omega (\tau _{1},\tau_2,\epsilon),
\end{align*}
is injective. Furthermore, the differential operators $dx^2 \boldsymbol{D}_{x}$ and $\nabla_x$ are equivariant in the sense that 
\begin{align}
dx^2 \boldsymbol{D}_{x}=
\left(F^{\Omega}\right)^{-1} \circ \nabla_x \circ F^{\Omega},
\quad 
\nabla_x\vert_{F^{\Omega}(\mathcal{D}_{\mathrm{sew}})}=F^{\Omega} \circ dx^2 \boldsymbol{D}_{x} \circ  \left(F^{\Omega}\right)^{-1}.
\label{eq:Dxnablaequiv}
\end{align}
\end{theorem}
\begin{remark}
\label{rem:DelD}
We will  write $dx^2 \boldsymbol{D}_{x}=\nabla_x$ below as shorthand for  \eqref{eq:Dxnablaequiv}.
\end{remark}

\begin{theorem}\label{thm:Dxnabla2}
The holomorphic map $
F^{\Omega }:\mathcal{D}_{\mathrm{sew}} \rightarrow \HH_{2}$ is not surjective.
\end{theorem}
\noindent \textbf{Proof.}
Assume that $F^{\Omega }$ is surjective and find a contradiction.
Consider the rank two Heisenberg VOA $M^2$ with partition function  (from \eqref{eq:Z2M})
\begin{align*}
Z_{M^2}^{(2)}(\tau _{1},\tau _{2},\epsilon ) = \frac{1}{\eta(\tau_1)^2\eta(\tau_2)^2\det\left( \mathbbm{1} - A_{1}A_2 \right)}.
\end{align*}
Since $F^{\Omega }$ is assumed surjective it follows that
$G_{1}(\Omega):=1/Z^{(2)}_{M^2}((F^{\Omega })^{-1}\Omega)$  is a solution on $ \HH_{2}$ to the differential equation (using  \eqref{eq:HeisDE})
\begin{align}
\left(\nabla_x + \frac{1}{6} s(x)\right) G_{1}(\Omega)=0.
\label{eq:ZMdiff}
\end{align}
For $\gamma\in \Sp(4,\Z)$, consider the modular transformation $\gamma:\Omega\rightarrow\widetilde{\Omega}$ of 
\eqref{eq:ZMdiff}. Then Lemma~\ref{lem:modder} implies that 
$\widetilde{G}_1(\Omega):=G_{1}(\widetilde{\Omega})/\det(C\Omega+D)$ is also a solution to 
\eqref{eq:ZMdiff}. Let $\chi= \widetilde{G}_1(\Omega)/G_1(\Omega)$ which must therefore  satisfy
$\nabla_x \chi=0$ since the $\{\Phi_r(x)\}$ of \eqref{eq:Phibasis} are independent holomorphic 2--differentials. Hence $\chi=\chi(\gamma)$ and it follows that
\begin{align*}
G_1(\widetilde{\Omega})=\chi(\gamma) \det(C\Omega + D) G_{1}(\Omega),
\end{align*}
i.e. 
$G_{1}(\Omega)$ is a meromorphic $\Sp(4,\Z)$ Siegel modular form of weight $1$ with a multiplier system $\chi(\gamma)$, a 1--dimensional complex character for $\Sp(4,\Z)$.  
The commutator subgroup of $\Sp(4,\Z)$ is of index $2$ in $\Sp(4,\Z)$
so that $\chi(\gamma)\in\{\pm 1\}$ for all $\gamma\in\Sp(4,\Z)$ \cite{Fr,Kl}. But
Theorem~\ref{Theorem_period_eps} implies that the left torus modular transformation $\tau_1\rightarrow \tau_1+1$ is equivalent to the $\Sp(4,\Z)$  transformation $T_1:\Omega_{11}\rightarrow \Omega_{11}+1$ with multiplier $\chi(T_1)=e^{i\pi/6}$ which contradicts that $\chi(T_1)\in\{\pm 1\}$. 
Hence $F^{\Omega }$  is not surjective. \hfill \qed

\begin{remark}
\label{rem:} We note that \eqref{eq:ZMdiff} is invariant under the $\Sp(4,\Z)$ subgroup  $\Gamma\simeq (SL(2,\Z)\times SL(2,\Z))\rtimes \Z_{2}$ of Theorem~\ref{Theorem_period_eps} which preserves $\mathcal{D}_{\mathrm{sew}}$ for  $G_{1}(\Omega)$ a weight $1$ Siegel form with a multiplier system $\chi(\gamma)\in\langle e^{i\pi/6}\rangle$ as shown in Theorem~8 of \cite{MT3}.
\end{remark}

\begin{remark}
We emphasise that $Z_{M^2}^{(2)}(\tau _{1},\tau _{2},\epsilon )$ is not a function on the full Siegel upper half plane but only on the image  $F^{\Omega }\left(\mathcal{D}_{\mathrm{sew}} \right)$.
\end{remark}

\subsection{A differential equation for holomorphic  $1$--differentials}
Consider two Heisenberg VOA  modules $M_{\lambda_1},M_{\lambda_2}$. 
The genus two $3$-point function for $h$ inserted at $x_1,x_2,y$ is  from Proposition~\ref{HeisenbergMnpoint} (and  \cite{MT3}) given by
\begin{align*}
& dx_1 dx_2 dy Z^{(2)}_{\boldsymbol{\lambda}} (h,x_1;h,x_2;h,y) \\
& = \Big( \nu_{\boldsymbol{\lambda}}(x_1)\nu_{\boldsymbol{\lambda}}(x_2)\nu_{\boldsymbol{\lambda}}(y) + \nu_{\boldsymbol{\lambda}}(x_1)\omega(x_2,y) \\ 
& \qquad + \nu_{\boldsymbol{\lambda}}(x_2)\omega(x_1,y) + \nu_{\boldsymbol{\lambda}}(y)\omega(x_1,x_2) \Big)
Z^{(2)}_{\boldsymbol{\lambda}} .
\end{align*}
Since $\widetilde{\omega}=\frac{1}{2}h[-1]^2\vac$  we find
\begin{align*}
& dx^2 dy Z^{(2)}_{\boldsymbol{\lambda}} (\widetilde{\omega},x;h,y) \\
&= 
\lim_{x_i \rightarrow x} \frac{1}{2}  dx_1 dx_2\left( dy\, Z^{(2)}_{\boldsymbol{\lambda}} (h,x_1;h,x_2;h,y) - \frac{\nu_{\boldsymbol{\lambda}}(y) }{(x_1-x_2)^2} Z^{(2)}_{\boldsymbol{\lambda}} \right) \\
&= \left( \frac{1}{2}\nu_{\boldsymbol{\lambda}}(x)^2\nu_{\boldsymbol{\lambda}}(y) + \nu_{\boldsymbol{\lambda}}(x)\omega(x,y) + \frac{1}{12} \nu_{\boldsymbol{\lambda}}(y)s(x) \right) Z^{(2)}_{\boldsymbol{\lambda}} .
\end{align*}
By Proposition~\ref{prop:WardId} and Theorem~\ref{thm:Dxnabla1} we also have 
\begin{align*}
& dx^2 dy Z^{(2)}_{\boldsymbol{\lambda}} (\widetilde{\omega},x;h,y) \\
&=
\Big( \nabla_x +  dx^2 \left( \sups{2}\mathcal{P}_{1}(x,y)\partial_y + \sups{2}\mathcal{P}_{2}(x,y) \right) \Big) Z^{(2)}_{\boldsymbol{\lambda}} (h,y) dy \\
&= 
\Big( \nabla_x +  dx^2 \left( \sups{2}\mathcal{P}_{1}(x,y)\partial_y + \sups{2}\mathcal{P}_{2}(x,y) \right) \Big)  \nu_{\boldsymbol{\lambda}}(y) Z^{(2)}_{\boldsymbol{\lambda}} .
\end{align*}
\noindent
Using \eqref{eq:Zlamom} and comparing we thus obtain (recalling that $\mathcal{P}_{2}(x,y)=\partial_y\mathcal{P}_{1}(x,y)$):
\begin{proposition}\label{prop:nablanu}
The genus two holomorphic 1--forms  $\nu_i(x)$, $i=1,2$,  satisfy the following differential equation on the sewing domain  $\mathcal{D}_{\mathrm{sew}}$
\begin{align}
\nabla_x \nu_i(y)+  dx^2 \partial_y \left( \sups{2}\mathcal{P}_{1}(x,y) \nu_i(y)\right)
=
\omega(x,y) \nu_i(x). 
\label{eq:nuDE}
\end{align}
\end{proposition}

Proposition~\ref{prop:nablanu}  allows us to determine a global analytic expression for the generalised Weierstrass function $\sups{2}\mathcal{P}_{1}(x,y)$:
\begin{proposition}\label{prop:P21}
$\sups{2}\mathcal{P}_{1}(x,y)$ is given by the $(2,-1)$--bidifferential
\begin{align}
\sups{2}\mathcal{P}_{1}(x,y)dx^2(dy)^{-1}=
-\frac{
\omega(x,y)\begin{vmatrix}
\nu_1(x) &\nu_1(y) \\
 \nu_2(x)& \nu_2(y)
\end{vmatrix}
+\begin{vmatrix}
\nu_1(y) &  \nabla_x\nu_1(y)\\
\nu_2(y) &  \nabla_x\nu_2(y)
\end{vmatrix}
}
{
\begin{vmatrix}
\nu_1(y) &  \partial_{y} \nu_1(y)\\
\nu_2(y) &  \partial_{y} \nu_2(y)
\end{vmatrix}
dy
},
\label{eq:P21form}
\end{align}
which is holomorphic for $x\neq y$ 
where, for any local coordinates $x,y$
\begin{align*}
\sups{2}\mathcal{P}_{1}(x,y)= \frac{1}{x-y}+\text{regular terms}.
\end{align*}
\end{proposition}
\noindent
\textbf{Proof.} 
Proposition~\ref{prop:nablanu} implies
\begin{align*}
  \nu_2(y)\nabla_x \nu_1(y)+  \nu_2(y)dx^2  \partial_y \left( \sups{2}\mathcal{P}_{1}(x,y)   \nu_1(y)\right)
&=
\nu_2(y)\omega(x,y) \nu_1(x),\\ 
\nu_1(y)\nabla_x \nu_2(y)+  \nu_1(y)dx^2  \partial_y \left( \sups{2}\mathcal{P}_{1}(x,y)   \nu_2(y)\right) 
&= 
\nu_1(y)\omega(x,y) \nu_2(x).
\end{align*}
Taking the difference we obtain \eqref{eq:P21form}. 
Thus $\sups{2}\mathcal{P}_{1}(x,y)dx^2(dy)^{-1}$ is  globally defined for all $\Omega\in\HH_2$  (i.e. not just on $F^{\Omega }(\mathcal{D}_{\mathrm{sew}})$).

Let $W(y)=\left|\begin{smallmatrix}
\nu_1(y) &  \partial_{y} \nu_1(y)\\
\nu_2(y) &  \partial_{y} \nu_2(y)
\end{smallmatrix}\right |dy $
denote the Wronskian denominator of  the right hand side of \eqref{eq:P21form}.   
$W(y)$ is a holomorphic $3$--differential  with $6$ zeros  in $y$ (counting multiplicity) from the  Riemann-Roch theorem \cite{FK}. 
The numerator of \eqref{eq:P21form} is a holomorphic  $(2,2)$--bidifferential  for $x\neq y$ with a simple pole at $x=y$ with residue  $-W(y)dx^2dy^{-1}$ so that 
$\sups{2}\mathcal{P}_{1}(x,y)\sim \frac{1}{x-y}$ for $x\sim y$.
But from \eqref{eq:nuDE},  $\mathcal{P}_{1}(x,y)$ cannot have any $x$--independent poles in $y$  so that the numerator of the right hand side of \eqref{eq:P21form} possesses the same $6$ zeros  in $y$ as $W(y)$. Thus $\sups{2}\mathcal{P}_{1}(x,y)$ is holomorphic for all $x\neq y$. \hfill \qed 
\begin{remark}\label{rem:Pijhol}
Proposition~\ref{prop:P21} implies that 
$\sups{2}\mathcal{P}_{i,1+j}(x,y)$ is similarly holomorphic for $x\neq y$ following \eqref{eq:P2i1j} with 
\begin{align*}
\sups{2}\mathcal{P}_{i,1+j}(x,y)=\frac{1}{(x-y)^{1+i+j}}+\text{regular terms}.
\end{align*}
\end{remark}
Referring to the genus two Zhu reduction of Theorem~\ref{thm:npZhured} and  combining Theorem~\ref{thm:Psi2forms}, 
Proposition~\ref{prop:P21} and Remark~\ref{rem:Pijhol} we conclude
\begin{theorem}
\label{th:Zhured2}
All the coefficients $\supsin{2}\mathcal{F}_a(x)$, $\supsin{2}\mathcal{F}^\Pi(x;1)$ and 
 $\sups{2}\mathcal{P}_{1+j}(x,y)$ involved in the $N=2$ genus two Zhu reduction are holomorphic for all $x\in \widehat{\mathcal{S}}_a$, $y\in \widehat{\mathcal{S}}_b$ for  $x\neq y$, and for all $(\tau_1,\tau_2,\epsilon)\in \mathcal{D}_{\mathrm{sew}}$. 
\end{theorem}
\begin{remark}
\label{rem:ZhuVir}
Theorem~\ref{th:Zhured2} implies that all coefficients appearing in the genus two  Ward identities of Propositions~\ref{prop:WardId} and \ref{prop:Vnpt} are convergent on $\mathcal{D}_{\mathrm{sew}}$. 
In particular, this implies that any genus two differential equation derived from a Virasoro singular vector has coefficients convergent on $\mathcal{D}_{\mathrm{sew}}$. 
This is explored further in \cite{GT1} and \cite{GT2}.
\end{remark}

\subsection{A differential equation for the normalised bidifferential}
Consider the genus two $4$-point function for $h$ inserted at $x_1,x_2,y_1,y_2 $ which from Proposition~\ref{HeisenbergMnpoint} (and  \cite{MT3}) is given by
\begin{align*}
& dx_1 dx_2 dy_1 dy_2 Z^{(2)}_M (h,x_1;h,x_2;h,y_1;h,y_2) \\
&= \Big( \omega(x_1,x_2)\omega(y_1,y_2) + \omega(x_1,y_1)\omega(x_2,y_2) \\
&\qquad + \omega(x_1,y_2)\omega(x_2,y_1) \Big) Z^{(2)}_M .
\end{align*}
Much as before, we find
\begin{align*}
& dx^2 dy_1 dy_2 Z^{(2)}_M (\widetilde{\omega},x;h,y_1;h,y_2) \\
&= \lim_{x_i \rightarrow x} \frac{1}{2} \bigg( dx_1 dx_2 dy_1 dy_2 Z^{(2)}_M (h,x_1;h,x_2;h,y_1;h,y_2) 
- \frac{dx_1 dx_2}{(x_1-x_2)^2} \omega(y_1,y_2) Z^{(2)}_M \bigg) \\
&= \left( \frac{1}{12}s(x) \omega(y_1,y_2) + \omega(x,y_1)\omega(x,y_2) \right) Z^{(2)}_M .
\end{align*}
By Proposition \ref{prop:WardId} and Theorem~\ref{thm:Dxnabla1}  we find
\begin{align*}
& dx^2 dy_1 dy_2 Z^{(2)}_M (\widetilde{\omega},x;h,y_1;h,y_2) \\
&= \Big( \nabla_x +dx^2  \sum_{r=1}^2 \left( \sups{2}\mathcal{P}_{1}(x,y_r)\partial_{y_r} + \sups{2}\mathcal{P}_{2}(x,y_r) \right) \Big) Z^{(2)}_M (h,y_1;h,y_2) dy_1 dy_2 \\
&= \Big( \nabla_x + dx^2 \sum_{r=1}^2 \left( \sups{2}\mathcal{P}_{1}(x,y_r)\partial_{y_r} + \sups{2}\mathcal{P}_{2}(x,y_r) \right) \Big) \omega(y_1,y_2) Z^{(2)}_M ,
\end{align*}
which implies on using Proposition~\ref{prop:HeisDE}
\begin{proposition}\label{nablaomega}
The bidifferential  $\omega(x,y)$ satisfies the differential equation 
\begin{align*}
\Big( \nabla_x + dx^2 \sum_{r=1}^2 \left( \sups{2}\mathcal{P}_{1}(x,y_r)\partial_{y_r} + \sups{2}\mathcal{P}_{2}(x,y_r) \right) \Big) \omega(y_1,y_2)
&=  \omega(x,y_1)\omega(x,y_2) .
\end{align*}
\end{proposition}
\noindent
This differential equation is similar in form to the genus one case \cite{HT}
\begin{align*}
&\Big( q\partial_q +
 \sum_{r=1}^2\left(P_1(x-y_r,\tau)\partial_{y_r}+ P_2(x-y_r,\tau) \right) \Big) P_2(y_1-y_2,\tau) \\
&= P_2(x-y_1,\tau)P_2(x-y_2,\tau).
\end{align*}

\subsection{A differential equation for the projective connection}
Very much as in the last example, we find
\begin{align*}
& dx^2 dy^2 Z^{(2)}_M (\widetilde{\omega},x;\widetilde{\omega},y) 
= \left( \frac{1}{144}s(x)s(y) +\frac{1}{2} \omega(x,y)^2 \right) Z^{(2)}_M .
\end{align*}
By Proposition \ref{prop:Vnpt} and Theorem~\ref{thm:Dxnabla1}  we find
\begin{align*}
& dx^2 Z^{(2)}_M (\widetilde{\omega},x;\widetilde{\omega},y) \\
&= \Big( \nabla_x +dx^2  \left( \sups{2}\mathcal{P}_{1}(x,y)\partial_{y} + \sups{2}\mathcal{P}_{2}(x,y) \right) \Big) Z^{(2)}_M (\widetilde{\omega},y)  
+\frac{1}{2}\cdot \sups{2}\mathcal{P}_{4}(x,y)  Z^{(2)}_Mdx^2.
\end{align*}
which implies on using Proposition~\ref{prop:HeisDE} that
\begin{proposition}\label{nablas(x)}
The genus two projective connection satisfies  
\begin{align*}
\Big(\nabla_x +  dx^2 \left(\sups{2}\mathcal{P}_{1}(x,y) \partial_y +2\sups{2}\mathcal{P}_{2}(x,y) \right)\Big) \left( \frac{1}{6}s(y)\right)+ \sups{2}\mathcal{P}_{4}(x,y)dx^2dy^2
=
\omega(x,y) ^2. 
\end{align*}
\end{proposition}

Finally, we note that Proposition~\ref{nablaomega} implies Proposition~\ref{prop:nablanu} 
on  integrating $y_1$ or $y_2$ over a $\beta^i$ cycle and applying \eqref{nui} and Proposition~\ref{nablas(x)} on taking the $y_1\rightarrow y$ limit.


\subsection{Conjectures}
We conclude with number of conjectures that naturally arise:
\begin{conjecture}\label{Nforms}
We conjecture that
\begin{enumerate}[(i)]
\item 
the formal inverse matrix $
\left( \mathbbm{1} - \widetilde{\Lambda}_a \widetilde{\Lambda}_{\bar{a}} \right)^{-1}
$
of \eqref{eq:IminLamLaminv}
is convergent on $\mathcal{D}_{\mathrm{sew}}$ for all $N\ge 2$,
\item
$\supsin{N}\mathcal{F}_a(x)dx^N$ and $\supsin{N}\mathcal{F}^\Pi(x;m)dx^N$ of 
Definition~\ref{calFforms}
form a dimension $2N-1$ basis of holomorphic $N$--differentials on $\mathcal{D}_{\mathrm{sew}}$  for all $N\ge 3$ (c.f. Remark~\ref{rem:ZhuTh3}),
\item
the generalised Weierstrass function$\sups{N}\mathcal{P}_{1}(x,y)$ 
of Definition~\ref{calP21} is holomorphic for all $x\in\widehat{\mathcal{S}}_a,y\in\widehat{\mathcal{S}}_{b}$ for  $x\neq y$ and $N\geq3$,
\item
$\sups{N}\mathcal{P}_{1}(x,y)dx^{N}(dy)^{1-N}$ is a globally defined  holomorphic $(N,1-N)$--bidifferential for $x\neq y$  for all $N\ge 3$.
\end{enumerate}
\end{conjecture}
\begin{conjecture}
The genus two partition function for a $C_2$-cofinite VOA is convergent on $\mathcal{D}_{\mathrm{sew}}$. In particular, a $C_2$-cofinite VOA partition function obeys a partial differential equation resulting \cite{GT1} from  genus two Zhu recursion applied to the genus two 1--point function for  a singular Virasoro vacuum descendant. This is explored in \cite{GT2} for the $c=-22/5$ Virasoro $(2, 5)$--minimal model.
\end{conjecture}
\begin{conjecture}
Proposition~\ref{prop:HeisDE} generalises to all genera in some particular sewing domain $\mathcal{D}^{(g)}_{\mathrm{sew}}$ where the genus $g$ partition function for the rank $1$ Heisenberg  VOA satisfies the differential equation
\begin{align*}
\left(\nabla_x -\frac{1}{12} s(x) \right)Z^{(g)}_M=0.
\label{eq:HeisDEg}
\end{align*}
for $x\in \mathcal{D}^{(g)}_{\mathrm{sew}}$ and 
where $\nabla_x$ is a suitable generalisation of \eqref{eq:nabla}, that is, a covariant derivative with respect to the surface moduli depending on a basis of holomorphic 2-differentials.
\end{conjecture}


\end{document}